\documentstyle{amsppt}

\magnification=\magstep1
\baselineskip=14pt
\parindent=0pt
\parskip=14pt

\vsize=7.7in
\voffset=-.4in
\hsize=5.7in

\overfullrule=0pt

\define\A{{\Bbb A}}
\define\C{{\Bbb C}}

\define\R{{\Bbb R}}
\define\Q{{\Bbb Q}}

\define\Z{{\Bbb Z}}


\redefine\H{\frak H}


\define\z{\frak z}

\define\a{\alpha}
\redefine\b{\beta}
\redefine\d{\delta}
\define\e{\epsilon}
\redefine\l{\lambda}
\redefine\o{\omega}
\define\ph{\varphi}

\define\s{\sigma}
\redefine\P{\Phi}
\predefine\Sec{\S}

\redefine\S{{\Cal S}}


\define\sh{\sharp}


\define\back{\backslash}

\define\lra{\longrightarrow}
\redefine\tt{\otimes}
\define\scr{\scriptstyle} 
\define\liminv#1{\underset{\underset{#1}\to\leftarrow}\to\lim}
\define\limdir#1{\underset{\underset{#1}\to\rightarrow}\to\lim}

\define\nass{\noalign{\smallskip}}


\define\CH{\widehat{CH}}


\define\db{\bar\partial}
\redefine\d{\partial}

\redefine\ord{\text{\rm ord}}
\define\Ei{\text{\rm Ei}}
\redefine\O{\Omega}
\predefine\oldvol{\vol}
\redefine\vol{\text{\rm vol}}

\define\pr{\text{\rm pr}}



\redefine\div{\text{\rm div}}

\define\Spec{\text{\rm Spec}}

\define\tr{\text{\rm tr}}

\define\sig{\text{\rm sig}}

%

%

\redefine\Im{\text{\rm Im}}
\redefine\Re{\text{\rm Re}}

\define\II{\int^\bullet}     

\define\HZ{\hat{\Z}}

\define\vth{\vartheta}

\define\GSpin{\text{\rm GSpin}}
\define\ev{C}

\define\CT#1{\operatornamewithlimits{CT}_{#1}}

\define\und#1{\underline{#1}}
\define\SL{\text{\rm SL}}
\define\sgn{\text{\rm sgn}}

\define\XX{\frak X}
\define\bXX{\bar\XX}
\define\ZZ{\frak Z}
\define\degh{\widehat{\text{\rm deg}}}
\define\Pich{\widehat{\text{\rm Pic}}}
\define\divh{\widehat{\text{\rm div}}}


\define\borchold{\bf1}
\define\borch{\bf2}
\define\borchduke{\bf3}
\define\borchdukeII{\bf4}
\define\bost{\bf5}
\define\bostgilletsoule{\bf6}
\define\bruinier{\bf7}
\define\bruinierII{\bf8} 
\define\bruinierkuehn{\bf9}
\define\burgoskuehn{\bf10}
\define\cohen{\bf11}
\define\EZ{\bf12}
\define\flenstedj{\bf13}
\define\freitaghermann{\bf13}
\define\funkethesis{\bf14}
\define\funke{\bf15}
\define\gelbart{\bf16}
\define\gilletsoule{\bf17}
\define\goldstein{\bf18}
\define\GN{\bf19} 
\define\harris{\bf20}
\define\harveymoore{\bf21}
\define\hejhal{\bf22}
\define\hermanna{\bf23}
\define\hermannb{\bf24}
\define\hirzebruchzagier{\bf25}
\define\duke{\bf26}
\define\annals{\bf27}
\define\kbourb{\bf28}
\define\kmillsonI{\bf28}
\define\kmillsonII{\bf29}
\define\kmcana{\bf30}
\define\krannals{\bf31}  
\define\kryII{\bf32} 
\define\ky{\bf33}
\define\kuehnthesis{\bf34}
\define\kuehn{\bf35}
\define\lebedev{\bf36}
\define\maillotroessler{\bf37}
\define\milne{\bf38}
\define\niebur{\bf39}
\define\petersson{\bf40}
\define\rademacher{\bf41}
\define\rademacherZ{\bf42}
\define\rallisHDC{\bf43}
\define\rallisschiff{\bf44}
\define\rohrlich{\bf45}
\define\siegeltata{\bf46}
\define\vdgeer{\bf47}
\define\vdgeerbook{\bf48}
\define\waldspurger{\bf49}
\define\weilI{\bf50}
\define\weilII{\bf51}
\define\zagier{\bf52}
\define\zuckerman{\bf53}

\centerline{\bf Integrals of Borcherds forms}
\smallskip
\centerline{\bf by}
\smallskip
\centerline{\bf Stephen S. Kudla\footnote{Partially 
supported by NSF grant DMS-9970506 and by a Max-Planck Research Prize 
from the Max-Planck Society and Alexander von Humboldt Stiftung. }}
\vskip .5in

\subheading{Introduction}

Let $V$ be a non-degenerate inner product space over $\Q$ of 
signature $(n,2)$, and let $D$ be the space of oriented negative 
$2$-planes in $V(\R)$. In \cite{\borch}, Borcherds
constructed certain meromorphic modular forms $\Psi(F)$ on $D$ with respect to 
arithmetic subgroups $\Gamma_M$ of $G=O(V)$ by regularizing the theta integral of 
vector valued elliptic modular forms $f$ of weight $1-\frac{n}2$ for $SL_2(\Z)$ 
with poles at the cusp, cf. also \cite{\borchold}, \cite{\harveymoore}, \cite{\bruinier}, \cite{\bruinierII}. 
The Borcherds forms $\Psi(f)$ can be viewed as 
meromorphic sections of powers of a certain line bundle $\Cal L$ on $X=\Gamma_M\back D$. 
Taking the standard Petersson metric $||\ ||$ on $\Cal L$, it is of interest 
in Arakelov geometry to compute the integral:
$$\kappa(\Psi(f)) := -\vol(X)^{-1}\int_{\Gamma_M\back D}\log||\Psi(z,f)||^2\,d\mu(z),\tag0.1$$
where $d\mu(z)$ is a $G(\R)$--invariant volume form on $D$.
 
In this paper, we give an explicit formula for $\kappa(\Psi(f))$ is almost all cases\footnote{
In the {\it exceptional cases}, where $\dim(V)=3$ (resp. $4$) and $V$ has Witt index $1$ (resp. $2$), 
some additional regularization is required.}. To describe it, suppose that $M$ is a lattice 
in $V$ such that the quadratic form $Q(x) = \frac12(x,x)$ is $\Z$--valued and 
let $M^\sh\supset M$ be the dual lattice. Recall that the modular form
$f$ is valued in the space $\C[M^\sh/M]$, for a suitable choice of $M$, and has a Fourier expansion of the form
$$f(\tau) = \sum_{\mu\in M^\sh/M} \sum_{m\in\Q} c_\mu(m)\,q^m\, \ph_\mu\tag0.2$$
where $\tau\in\H$, $q^m=e(m\tau)$, and where $c_\mu(m)$ is zero unless $m\in Q(\mu)+\Z$ 
and $m>-R$ for some positive integer $R$. In addition, if $m\le0$, then $c_\mu(m)\in\Z$.
Let 
$$\Gamma_M=\{ \gamma\in SO(V)(\Q)\mid \gamma M = M \text{ and $\gamma$ acts trivially in $M^\sh/M$}\},\tag0.3$$
and let $X_M = \Gamma_M\back D$, so that $X$ is a quasi--projective variety. 
For each $m>0$ and $\mu\in M^\sh/M$, there is a divisor $Z(m,\mu)$ on $X$, associated to 
the set of vectors $x\in \mu+M$ with $Q(x)=m$. These divisors, called rational quadratic divisors
or Heegner divisors in \cite{\borch}, include the Heegner points, for $n=1$, the
Hirzebruch--Zagier curves on Hilbert modular surfaces, for $n=2$, and the Humbert surfaces on 
Siegel threefolds, for $n=3$. They are also special cases of the cycles considered in
\cite{\kmillsonI}, \cite{\kmillsonII}, \cite{\duke}, etc.. 
Then, a key fact, due to Borcherds \cite{\borch}, is that the divisor of the form $\Psi(f)^2$, 
which has weight $c_0(0)$, is an
explicit linear combination of these cycles:
$$\div(\Psi(f)^2) = \sum_{\mu}\sum_{m>0} c_\mu(-m)\,Z(m,\mu).\tag0.4$$

Our first result concerns the generating function for the degrees of the cycles $Z(m,\mu)$. 
Let $\O$ be the first Chern form of the metrized line bundle $\Cal L^\vee$ 
on $X$, dual to $\Cal L$, and let
$$\deg(Z(m,\mu)) = \int_{Z(m,\mu)} \O^{n-1}\tag0.5$$
be the volume of the cycle $Z(m,\mu)$ with respect to $\O$. Similarly, 
let
$$\vol(X) = \int_X \O^n.\tag0.6$$

From now on, we exclude the two exceptional cases. In addition, for simplicity here in the introduction, 
we make the following
`class number $1$' assumption\footnote{Without this assumption, one must either work adelically, as is done in the body of the 
paper, or introduce sums over a suitable collection of lattices determined by the 
decomposition (1.3)}: Let
$$K_M=\{g\in SO(V)(\A_f)\mid gM=M\ \text{ and $g$ acts trivially on $M^\sh/M$}\ \}$$
and suppose that
$$SO(V)(\A_f) = SO(V)(\Q)^+\,K_M,$$
where $SO(V)(\Q)^+ = SO(V)(\Q)\cap SO(V)(\R)^+$, for $SO(V)(\R)^+$ the identity 
component of $SO(V)(\R)$.

We then show the following, using the Siegel--Weil formula and results of 
\cite{\kmillsonI}, \cite{\kmillsonII}, \cite{\kmcana}.
\proclaim{Theorem I} For each $\mu\in M^\sh/M$, there 
is an Eisenstein series $E(\tau,s;\mu,\frac{n}2+1)$, for $\tau\in \H$ and $s\in \C$, of weight $\frac{n}2+1$ 
such that 
$$E(\tau,s_0;\mu) = \delta_{\mu,0}+ \vol(X)^{-1}\sum_{m>0} \deg(Z(m,\mu))\,q^m,$$
where $s_0=\frac{n}2$. 
\endproclaim

Analogous results for generating functions for cycles of higher codimension for 
anisotropic $V$'s are discussed in \cite{\kbourb}, section 3, and in \cite{\duke}. 

The integral $\kappa(\Psi(F))$ can also be expressed using the second term in the Laurent expansion 
at $s_0=\frac{n}2$ of these Eisenstein series. 
\proclaim{Theorem II} For each $\mu\in M^\sh/M$, the Fourier coefficients in the 
expansion
$$E(\tau,s;\mu) = \sum_m A_\mu(s,m,v)\,q^m$$
have Laurent expansion at $s=s_0 = \frac{n}2$
$$A_\mu(s,m,v) = a_\mu(m) + b_\mu(m,v)(s-s_0) + O((s-s_0)^2).$$
Let
$$\kappa_\mu(m) = \cases \lim\limits_{v\rightarrow\infty} b_\mu(m,v)&\text{ if $m>0$,}\\
\nass
\frac12\,(\log(2\pi)-\gamma)&\text{ if $m=0$,}\\
\nass
0&\text{ if $m<0$.}
\endcases
$$
Then, for $f$ with Fourier expansion (0.2), 
$$\kappa(\Psi(f)) = \sum_{\mu}\sum_{m\ge0} c_\mu(-m)\,\kappa_\mu(m).$$
\endproclaim
In addition, we derive the useful relation
$$-\vol(X)\,c_0(0) = \sum_{\mu}\sum_{m>0} c_\mu(-m)\,\deg(Z(m,\mu)).$$
The quantities $\kappa_\mu(m)$ can be calculated quite readily in any particular case; this will be done 
in a sequel \cite{\ky}. 

As an illustration, consider the case where 
$M=\Z^5$ with quadratic form of signature $(3,2)$ 
defined by $Q(x) =\frac12 {}^tx Q x$ where
$$Q = \pmatrix {}&1&{}&{}&{}\\1&{}&{}&{}&{}\\{}&{}&2&{}&{}\\{}&{}&{}&{}&1\\{}&{}&{}&1&{}\endpmatrix.$$
In this example, which is worked out in detail in section 5, 
$|M^\sh/M|=2$ and, labeling the cosets by $\mu=0$, $1$, we have
$$E(\tau,\frac32;\mu) = \delta_{\mu,0} + 
\zeta(-3)^{-1}\sum_{m>0\atop 4m\equiv \mu\mod(4)}
H(2,4m)\,q^m,$$
where $H(2,N)$ is the $N$-th coefficient in Cohen's Eisenstein series of weight $\frac52$, \cite{\cohen},
$$\Cal H_2(\tau) = \zeta(-3) + \sum_{N>0\atop N \equiv 0,1\mod(4)} H(2,N)\,q^N.$$ 
In this case, as explained in \cite{\vdgeer} and \cite{\GN}, 
$\Gamma_M\back D\simeq \text{\rm Sp}_4(\Z)\back \H_2$ is the Siegel threefold of level $1$, 
$\vol(X) = \zeta(-1)\zeta(-3)$, 
and $Z(m,\mu)$, for $4m\equiv \mu\mod(4)$, is the Humbert surface $\Cal G_{4m}$, in the notation of 
\cite{\vdgeer}.
Thus, Theorem I implies that
$$\deg(H_N) = -\frac1{12}\,H(2,N),$$
a relation due to van der Geer, \cite{\vdgeer}. Also, we find that, 
for $m>0$ with $4m=n^2d$ for a fundamental discriminant $d$, and with $4m\equiv \mu\mod(4)$, 
$$\align
\kappa_\mu(m) &=  \zeta(-3)^{-1}\,H(2,4m)\,\bigg[\frac43 \, 
+2\,\frac{\zeta'(-3)}{\zeta(-3)} -\frac12\log(d)-\frac{L'(-1,\chi_d)}{L(-1,\chi_d)} -C\\
\nass
\nass
{}&\qquad\qquad\qquad\qquad\qquad\qquad +\sum_{p\mid n} \bigg(\,\log|n|_p-\frac{b'_p(n,-1)}{b_p(n,-1)}\, \bigg) \,\bigg]. 
\endalign
$$
If $4m\not\equiv \mu\mod(4)$, then $\kappa_\mu(m)=0$. Here $L(s,\chi_d)$ is the L-series 
for the quadratic character $\chi_d$ and the other quantities are explained in section 5. 
It is shown by Gritesenko and Nikulin \cite{\GN} that the Siegel cusp $\Delta_5$ of weight $5$ and quadratic character
arises as a Borcherds form $\Psi(\bold f_5) = 2^{-6}\,\Delta_5(z)$, for a suitable meromorphic form $\bold f_5$ of weight
$-\frac12$ with expansion
$$\bold f_5(\tau) = \big(\ 10 + 108\, q + 808\, q^2 + \dots\ \big)\,\ph_0 + \big(\ 
q^{-\frac14} - 64\, q^{\frac34} - 513\, q^{\frac74} +\dots\ \big)\,\ph_1.$$
Thus, by Theorem II,
$$\multline
-\vol(X)^{-1}\int_X \log\big(|\Delta_5(z)|^2\,\det(y)^{5}\big)\cdot\O^3\\
\nass
 = 10\,\bigg[ - \frac43
-2\,\frac{\zeta'(-3)}{\zeta(-3)}+\frac{\zeta'(-1)}{\zeta(-1)} +\frac32\,\log(2)+\log(\pi)\,\bigg] - 7\log(2).
\endmultline$$

The main idea in the proof of Theorem II is the following. Recall that, in Borcherds 
construction, it is essentially the quantity  
$\log||\Psi(f)||^2$, rather than the meromorphic form $\Psi(f)$ itself, which arises as a regularized theta 
integral. Therefore, after some justification, we can compute the integral of this quantity by 
first integrating the theta kernel over $X$ and then
taking the regularized integral against $f$. This procedure is valid provided the integral of the
theta kernel is termwise abosolutely convergent, and it is for this reason that the exceptional
cases must be excluded. The Siegel--Weil formula
then identifies the integral of the 
theta kernel as a special value of an Eisenstein series of weight $\frac{n}2-1$ at the point $s_0 = \frac{n}2$. The regularized 
integral of this series against $f$ can then be evaluated by using 
Maass operators, which shifts the weight to $\frac{n}2+1$, and a Stokes theorem argument from section 9 of \cite{\borch}.

In fact, the method used here should also be applicable to the calculation of the integrals of the 
functions arising via Borcherds's construction for more general signatures $(p,q)$, and it would be 
interesting to investigate such cases. Note, in particular, that the remarkable product formulas for the
$\Psi(F)$'s in the case of signature $(n,2)$ play no role.

Possible applications of the formula for $\kappa(\Psi(f))$ to arithmetic geometry are discussed 
in section 6. The main point is that there should be a close connection between the second term in 
the Laurent expansion of the Fourier coefficients of the Eisenstein series $E(\tau,s;\mu)$ at $s_0=\frac{n}2$, and the heights 
of the divisors $Z(m,\mu)$ on $X$, after extension to a suitable integral model. Such a connection is also suggested by 
the results of joint work \cite{\kryII} with Michael Rapoport and Tonghai Yang in which we compute the 
heights of Heegner type divisors on the arithmetic surfaces $\frak X$ defined by Shimura curves, the case
$n=1$ with $V$ anisotropic. 
In fact, for suitablly defined classes $\widehat{\frak Z}(m,v)\in \CH^1(\frak X)$, the arithmetic Chow group
\footnote{The definition of the class $\widehat{\frak Z}(m,v)$ is still somewhat provisional.}  
of $\frak X$, and for a normalized version $\Cal E(\tau,s;\ph)$ of the Eisenstein 
series $E(\tau,s;\ph)$ of weight $\frac32$, we show that
$$\Cal E'(\tau,\frac12;\ph) = \sum_{m} \langle\,\widehat{\frak Z}(m,v), \hat\o\,\rangle\,q^m,$$
where $\tau = u+iv$, $\hat\o\in\CH^1(\frak X)$ is an extension of the metrized line bundle $\Cal L^\vee$, dual to $\Cal L$ to $\frak X$, 
and $\langle\ ,\ \rangle$ is the Gillet--Soul\'e height pairing. Thus, the second term in the 
Eisenstein series gives a generating functions for the `arithmetic volumes', at least in this example.

Here is a summary of the contents of the present paper. 
In section 1, we review the construction of the Borcherds forms $\Psi(F)$. 
An adelic formulation of this construction is given, which allows us to work 
more easily for general lattices and to make use of the adelic formulation of the 
Siegel--Weil formula and representation theory. Some explanation is given about how to pass
back and forth between the adelic and classical version. In section 2, derive the 
formula for $\Psi(f)$, assuming certain facts about Eisenstein series, the Siegel--Weil 
formula, and about convergence. In section 3, we consider convergence questions and, 
in particular, justify the interchange of the integration of the theta kernel with the Borcherds regularized
integral. In section 4, we first review the case of the Siegel--Weil formula which we need, 
including a refinement, already described by Weil, which is crucial in relating the 
integral over the orthogonal group occuring in this formula with the geometric integral 
we actually encounter. We then describe a general matching principle and apply it, together 
with the theory of \cite{\kmillsonI}, \cite{\kmillsonII}, \cite{\kmcana}, to prove that 
the degree generating function is given by the value of our Eisenstein series of weight $\frac{n}2+1$, 
as in Theorem I. This matching principle implies the coincidence of theta integrals for 
different quadratic spaces. For example, it shows that the degrees of the cycles $Z(m,\mu)$ 
occurring for spaces of signature $(n,2)$ always 
coincide with certain weighted representation numbers for spaces of signature $(n+2,0)$. 
This principle should have many other interesting applications.
In section 5, we discuss the example of signature $(3,2)$ described above. 
In section 6, we give some speculations about the applications of the formulas for $\kappa(\Psi(f))$'s 
in arithmetic geometry. This section also contains a brief discussion of the relation 
between our formula and work of Rohrlich \cite{\rohrlich} on analogous integrals of elliptic modular forms.

Work on the possibility of using Borcherds' forms $\Psi(F)$ in Arakelov theory began at
the program on Arithmetic Geometry the Isaac Newton Institute during May--June, 1998. 
The main steps in computing $\kappa(\Psi(F))$ were done during a stay at Orsay in June of 1999. 
The examples in section 5 were worked out during a visit to Humbolt University 
in Berlin in June of 2001. The speculations in section 6 profited from discussions with Ulf 
K\"uhn at that time. I would like to thank these institutions and my hosts (J. Nekovar and C. Soul\'e in Cambridge, 
J.-B. Bost and G. Henniart in Orsay, and J. Kramer in Berlin) for providing a wonderful working environment. 

I would like to thank A. Abbes, R. Borcherds, J.-B. Bost, Jens Funke, M. Harris, J. Kramer, Ulf K\"uhn,  
J. Millson, J. Nekovar, M. Rapoport, D. Rohrlich, E. Ullmo and T. Yang for stimulating discussions and valuable
suggestions. 
I would particularly like to thank Tonghai Yang for allowing me to quote the results of our joint project on the 
derivatives of Fourier coefficients of Eisenstein series and for many incisive comments, which considerably 
improved this paper. 

Finally, this work has been supported by NSF grant DMS-9970506 and by a Max-Planck Research Prize 
from the Max-Planck Society and the Alexander von Humboldt-Stiftung. 

\subheading{Contents}
{\obeylines\parskip=5pt\bf
\Sec1. Borcherd's forms
\Sec2. Computation of a regularized integral
\Sec3. Convergence estimates
\Sec4. Formulas for degrees
\Sec5. Examples
\Sec6. Speculations
}

\subheading{\Sec1. Borcherds' forms}

In this section we give an adelic formulation of a result of Borcherds on the construction of 
meromorphic modular forms. This formulation is convenient from the point of view of 
Hecke operators and Shimura varieties. Moreover, it is essential if we want to make use of 
the adelic version of the Siegel--Weil formula. 

Let $V$ be a vector space over $\Q$ with a non-degenerate quadratic form of signature $(n,2)$, 
and let $H=\GSpin(V)$. We write $(x,y) = Q(x+y)-Q(x)-Q(y)$ 
for the associated bilinear form. Let $D$ be the space of oriented negative $2$--planes in $V(\R)$. 
Recall that $D$ is isomorphic to the open subset $Q_-$ of the quadric $Q\subset \Bbb P(V(\C))$ defined 
by 
$$Q_- = \{w\in V(\C)\mid (w,w) =0, \ (w,\bar w)<0\}/\C^\times.\tag1.1$$
The isomorphism is given by $z \mapsto v_1-iv_2 = w$, where $v_1$, $v_2$ is a properly oriented 
basis for $z\in D$ with $(v_1,v_1)=(v_2,v_2)=-1$ and $(v_1,v_2)=0$. 
For a compact open subgroup $K\subset H(\A_f)$,  the space
$$X_K = H(\Q)\back\bigg( D\times H(\A_f)/K\bigg)\tag1.2$$
is the set of complex points of a quasi-projective variety rational over $\Q$ (via canonical models). 
This variety is projective if and only if 
$V$ is anisotropic over $\Q$. It is smooth if the image of $K$ in $SO(V)(\A_f)$ is neat.
Fix a component $D^+$ of $D$, and write
$$H(\A) =\coprod_j H(\Q) H(\R)^+ h_j K,$$
where $H(\R)^+$ is the identity component of $H(\R)\simeq \GSpin(n,2)$. 
Then
$$X_K \simeq \coprod_j \Gamma_j\back D^+,\tag1.3$$
where $\Gamma_j=H(\Q)\cap \big(\ H(\R)^+h_j K h_j^{-1}\ \big)$.

Let $\Cal  L_D$ be the restriction to $D\simeq Q_-$ of the tautological bundle on $\Bbb P(V(\C))$. The
action of $O(V)(\R)$ on $V(\C)$ induces an action of $H(\R)^+$ on $\Cal L_D$, and hence there is a holomorphic 
line bundle
$$\Cal L = H(\Q)\back \bigg( \Cal L_D\times H(\A_f)/K\bigg) \lra X_K.\tag1.4$$
This line bundle is also algebraic and has a canonical model over $\Q$, \cite{\harris}. On the component $\Gamma_j\back D^+$, 
$\Cal L$ has the form $\Gamma_j\back \Cal L_D$. 
Define a Hermitian metric $h_{\Cal L}$ on $\Cal L_D$ by taking 
$$h_{\Cal L}(w_1,w_2) = -\frac12 (w_1,\bar w_2).\tag1.5$$
This metric is clearly invariant under the natural action of $O(V)(\R)$ and hence descends to $\Cal L$.  

For a Witt decomposition
$$V(\R) = V_0 + \R e+\R f,\tag1.6$$
where $e$ and $f$, with $(e,f)=1$ and $(e,e)=(f,f)=0$, span a hyperbolic plane 
with orthogonal complement $V_0$, note that 
$\sig(V_0)=(n-1,1)$ and let 
$$C=\{ y\in V_0\mid (y,y)<0\}\tag1.7$$ 
be the negative cone. Then $D\simeq Q_-$ is isomorphic to the tube domain 
$$\Bbb D =\{z\in V_0(\C)\mid y= \Im(z)\in C\},\tag1.8$$
via the map 
$$\Bbb D \lra V(\C),\qquad z \mapsto w(z) := z+e-Q(z)f.\tag1.9$$
composed with the projection to $Q_-$. 
The map $z\mapsto w(z)$ can be viewed as a nowhere vanishing holomorphic section of $\Cal L_D$. 
Note that this section has norm
$$||w(z)||^2 = -\frac12\,(w(z),\bar w(z)) = -(y,y) =: |y|^2.\tag1.10$$
For $h\in O(V(\R))$ or $H(\R)$, we have
$$h\cdot w(z) = w(hz)\, j(h,z)\tag1.11$$
for a holomorphic automorphy factor 
$$j: H(\R)\times D \lra \C^\times.\tag1.12$$
For $k\in \Z$, holomorphic sections of $\Cal L^{\tt k}$ can be identified with
holomorphic functions 
$$\Psi: D\times H(\A_f) \lra \C\tag1.13$$ 
such that $\Psi(z,hk) =\Psi(z,h)$ for all $k\in K$ 
and 
$$\Psi(\gamma z,\gamma h) = j(\gamma,z)^k\,\Psi(z,h)\tag1.14$$
for all $\gamma\in H(\Q)$. The norm of the section $(z,h)\to \Psi(z,h)\cdot w(z)^{\tt k}$ associated to $\Psi$ is 
then
$$||\Psi(z,h)||^2 = |\Psi(z,h)|^2 \,|y|^{2k}.\tag1.15$$ 
We will refer to this as the Petersson norm of $\Psi$. Note that, under the isomorphism (1.3), $\Psi$ 
corresponds to the collection $(\Psi(\cdot, h_j))_{\{j\}}$ of holomorphic functions on $D^+$ 
automorphic of weight $k$ with respect to the $\Gamma_j$'s.  

{\bf Remark:} In the case $n=1$, so that $\Bbb D = \H^+\cup\H^-$, the automorphy 
factor is
$$j(g,z)=\det(g)^{-1}(cz+d)^2,$$
so that the `classical weight' of a section of $\Cal L^{\tt k}$ is $2k$.

We now give a version of Borcherds' construction \cite{\borch} of meromorphic sections of (a certain twist of) 
$\Cal L^{\tt k}$. These are obtained by a regularized theta lift for the 
dual pair $(SL_2, O(V))$.  

The basic theta kernel is constructed as follows. Let $S(V(\A))$, $S(V(\A_f))$, and $S(V(\R))$ be the 
Schwartz spaces of $V(\A)$, $V(\A_f)$, and $V(\R)$ respectively.  For $z\in D$, let 
$\pr_z:V(\R)\rightarrow z$ be the projection with kernel $z^\perp$, 
and, for $x\in V(\R)$, let 
$$R(x,z) = -(\pr_z(x),\pr_z(x)) = |(x,w(z))|^2|y|^{-2}.\tag1.16$$ 
Then the majorant associated to $z$ is 
$$(x,x)_z = (x,x) + 2R(x,z),\tag1.17$$
and the Gaussian is the function
$$\ph_\infty\in S(V(\R))\tt A^0(D),\qquad \ph_\infty(x,z) = e^{-\pi(x,x)_{\scr z}}.\tag1.18$$
Here $A^0(D)$ is the space of smooth functions on $D$. Note that, for $h\in O(V(\R))$, 
$$\ph_\infty(hx,hz)=\ph_\infty(x,z).\tag1.19$$

Let $G=SL_2$ and let $G'_\A$ be the $2$-fold metaplectic cover of $G(\A)$. Let 
$G'_\Q\subset G'_\A$ be the image of $G(\Q)$ under the canonical splitting homomorphism. 
The group $G'_\A$ acts in $S(V(\A))$ via the Weil representation $\o$ (determined by the 
standard additive character $\psi$ of $\A/\Q$ such that $\psi_\infty(x)=e(x)=e^{2\pi i x}$\,) and this action commutes with the linear action of $O(V)(\A)$.
It will sometimes be convenient to write this linear action as 
$\o(h)\ph(x) = \ph(h^{-1}x)$. 
For $z\in D$, $h\in O(V)(\A_f)$ and $g'\in G'_\A$, we let $\theta(g',z,h)$ be the linear functional 
on $S(V(\A_f))$ defined by
$$\ph\mapsto \theta(g',z,h;\ph) = \sum_{x\in V(\Q)} \o(g')\bigg(\ph_{\infty}(\cdot,z)\tt\o(h)\ph\bigg)(x).\tag1.20$$
Then, for $\gamma\in O(V)(\Q)$, we have
$$\theta(g',\gamma z,\gamma h;\ph) = \theta(g',z,h;\ph).\tag1.21$$
Also, by Poisson summation, \cite{\weilI}, for $\gamma\in G'_\Q$,  
$$\theta(\gamma g',z,h;\ph) = \theta(g',z,h;\ph).\tag1.22$$
Finally, for $g_1'\in G'_{\A_f}$ and $h_1\in O(V)(\A_f)$, we have
$$\theta(g'g'_1,z,hh_1;\ph)=\theta(g',z,h;\o(g'_1,h_1)\ph).\tag1.23$$
In particular, if $K\subset H(\A_f)$ is as above, and if $\ph\in S(V(\A_f))^K$, 
then the function
$$(z,h) \mapsto \theta(g',z,h;\ph)\tag1.24$$
on $D\times H(\A_f)$ descends to a function on $X_K$. We may view it as a linear 
functional on the space $S(V(\A_f))^K$ and hence we obtain: 
$$\align
\theta: G'_\Q\back G'_\A\times X_K &\lra \bigg(S(V(\A_f))^K\bigg)^\vee.\tag1.25\\
\nass
(g',z,h)&\mapsto \theta(g',z,h;\cdot).
\endalign
$$
Note that this function is {\it not} holomorphic in $z$. 

Let $K'_\infty$ be the full inverse image of $SO(2)\subset SL_2(\R)=G(\R)$ in $G'_\R$, 
for each $r\in \frac12\Z$ let $\chi_r$ be the character of $K'_\infty$ 
such that
$$\chi_\ell(k')^2 = e^{2ir\theta},\qquad
\text{if}\quad k'\mapsto k_\theta = \pmatrix \cos(\theta)&\sin(\theta)\\-\sin(\theta)&\cos(\theta)\endpmatrix \in SO(2)\tag1.26$$
under the covering projection. 
Let $K'\subset G'_\A$ be the full inverse image of $SL_2(\hat\Z)\subset G(\A_f)$, 
and note that 
$$G'_\A=G'_\Q G'_\R K'.\tag1.27$$ 
The Gaussian (1.18) is an eigenfunction of $K'_\infty$ with 
$$\o(k'_\infty)\ph_\infty(x,z) = \chi_\ell(k'_\infty)\ph_\infty(x,z),\tag1.28$$
for $\ell = \frac{n}2-1$. 
It then follows from (1.23) that 
$$\theta(g'k'_\infty k',z,h) = \chi_\ell(k'_\infty)\,(\o(k')^\vee)^{-1}\theta(g',z,h)\tag1.29$$
for all $k'_\infty\in K'_\infty$ and $k'\in K'$. In particular, the theta function has weight $\ell =\frac{n}2-1$.
Here $\o(k')^\vee$ denotes the action of $K'$ on the space $S(V(\A_f))^\vee$ 
dual to its action on $S(V(\A_f))$.

Now suppose that $F:G'_\Q\back G'_\A\rightarrow S(V(\A_f))^K$ is a function such
that 
$$F(g'k'_\infty k') = \chi_{-\ell}(k'_\infty)\,\o(k')^{-1} F(g')\tag1.30$$
for all $k'_\infty\in K'_\infty$ and $k'\in K'$. Then, as a function of $g'$, the 
$\C$--bilinear pairing
$$(\!(\ F(g'), \theta(g',z,h)\ )\!) = \theta(g',z,h;F(g'))\tag1.31$$
is left $G'_\Q$--invariant and right $K'_\infty K'$--invariant. 
Its integral over $G'_\Q\back G'_\A$, defined in general by a suitable regularization, is a function 
$$\P(z,h;F) = \II_{G'_\Q\back G'_\A} (\!(\ F(g'), \theta(g',z,h)\ )\!)\, dg'\tag1.32$$
on $X_K$.  

The Borcherds forms \cite{\borch} arise when $F$ comes from a certain type of vector valued
automorphic form {\it with possible poles at the cusps}. To describe these, it is convenient 
to pass to a point of view intermediate between that just explained and the 
classical formulation. 

Observe that $G'_\Q\cap (G'_\R K') \simeq SL_2(\Z)$. 
Let $\Gamma'$ be the full inverse image of $SL_2(\Z)\subset SL_2(\R)=G(\R)$ in 
the metaplectic cover $G'_\R$. Thus $\Gamma'$ is an extension of $SL_2(\Z)$ by 
$\{\pm1\}$. 
For each $\gamma'\in \Gamma'$, with image $\gamma$ in $SL_2(\Z)$, there is a unique element $\gamma''$ 
such that $\gamma'\gamma'' =\gamma \in G'_\Q\cap (G'_\R K')$.
For $\tau=u+iv\in \frak H$, the upper halfplane, let
$$g_\tau = \pmatrix 1&u\\{}&1\endpmatrix \pmatrix v^{\frac12}&{}\\{}&v^{-\frac12}\endpmatrix,\tag1.33$$
and let $g'_\tau= [g_\tau,1]\in G'_\R$. 
We then have
$$\gamma' g'_\tau = g'_{\gamma(\tau)} k'_\infty(\gamma',\tau)\tag1.34$$
for a unique $k'_\infty(\gamma',\tau)\in K'_\infty$. 
For $r\in \frac12\Z$, define an automorphy factor by 
$$j_r:\Gamma'\times \frak H \rightarrow \C^\times, 
\qquad\quad j_r(\gamma',\tau) = \chi_{-r}(k'_\infty(\gamma',\tau))\, |c\tau+d|^r,\tag1.35$$
if $\gamma=\pmatrix a&b\\c&d\endpmatrix$. For example, if $r\in \Z$, 
$j_r(\gamma',\tau) = (c\tau+d)^r$.

\proclaim{Lemma 1.1}
Suppose that $(\rho,\Cal V)$ is a representation of $K'$ and that
$$\phi:G'_\Q\back G'_\A \lra \Cal V$$
is a function such that 
$$\phi(g'k'_\infty k') = \chi_r(k'_\infty)\,\rho(k')^{-1}\, \phi(g').$$
Let
$$f(\tau) = v^{-r/2}\,\phi(g'_\tau).$$
Then, for all $\gamma = \gamma'\gamma''\in SL_2(\Z)$, 
$$f(\gamma(\tau)) = j_r(\gamma',\tau)\,\rho(\gamma'')f(\tau).$$
\endproclaim
\demo{Proof} We have
$$\eqalign{
f(\gamma(\tau))&=  v(\gamma(\tau))^{-r/2}\,\phi( g'_{\gamma(\tau)})\cr
\nass
\nass
{}&=|c\tau+d|^r\, v^{-r/2}\,\phi(\gamma g'_\tau k'_\infty(\gamma',\tau)^{-1}\,(\gamma'')^{-1})\cr
\nass
\nass
{}&=|c\tau+d|^r\,\chi_{-r}(k'_\infty(\gamma',\tau))\, v^{-r/2}\,\rho(\gamma'')\,\phi(g'_\tau)\cr
\nass
\nass
{}&=j_r(\gamma',\tau)\,\rho(\gamma'')f(\tau),\cr}\tag1.36
$$
as claimed.
\enddemo

Note that we can view $\Cal V$ as a representation of $\Gamma'$ by setting $\rho(\gamma')=\rho(\gamma'')$. 

Applying Lemma~1.1, via (1.29) and (1.30), we obtain automorphic forms
$$\vth(\tau,z,h) = v^{-\ell/2}\theta(g'_\tau,z,h),\tag1.37$$
of weight $\ell$, and
$$f(\tau) = v^{\ell/2}F(g'_\tau),\tag1.38$$
of weight $-\ell$, valued in $S(V(\A_f))^\vee$ and $S(V(\A_f))^K$ respectively. 
Note that $\vth$ is not holomorphic in $\tau$. 
Then the quantity in (1.32) is given by 
$$\P(z,h;F) = \II_{SL_2(\Z)\back\frak H} (\!(\ f(\tau), \vth(\tau,z,h)\ )\!)\ v^{-2}\,du\,dv\tag1.39$$
for a suitable choice of measure on $G'_\Q\back G'_\A$.

Let $M$ be a $\Z$--lattice in $V$, on which the quadratic form $Q(x) =\frac12(x,x)$
takes integral values, and let $M^\sharp$ be the dual lattice. 
Let $S_M\subset S(V(\A_f))$ 
be the space of functions with support in $\hat M^\sharp:=M^\sh\tt_\Z\HZ$ and 
constant on cosets of $\hat M:=M\tt_\Z\HZ$. We will 
use the characteristic functions of cosets as a basis for this 
finite dimensional space. 
The space $S_M$ is stable under the action of $K'$. The restriction to $S_M$ of the 
theta function $\vth(\tau,z,h)$, viewed as a linear functional, defines a (non-holomorphic) modular form 
of weight $\ell=\frac{n}2-1$ valued in $(\o^\vee,S_M^\vee)$, the dual of the representation 
$(\o,S_M)$ of $K'$.

Suppose that $F$ and, hence, $f$ takes values in $S_M$ and is meromorphic at the cusp in the following 
sense. Write
$$f(\tau) = \sum_{\ph} f_\ph(\tau)\cdot \ph,\tag1.40$$
where $\ph$ runs over the coset basis for $S_M$, and let
$$f_\ph(\tau) =\sum_{m\in \Q} c_\ph(m)\, q^m\tag1.41$$
be the Fourier expansion of $f_\ph$, where $q^m=e(m\tau)$. 
We will sometimes write $c_0(m)$ for the Fourier coefficients of $f_{\ph_0}$ 
where $\ph_0$ is the characteristic function of $\hat M$; the constant term $c_0(0)$ will 
play a crucial role. 
The Fourier coefficients $c_\ph(m)$ are nonzero only for $m\in \frac1{N}\Z$, for some integer $N$, 
and we require that only a finite number of $c_\ph(m)$'s with $m<0$ are nonzero.  
Then the pairing
$$(\!(\ f(\tau), \vth(\tau,z,h)\ )\!) = \sum_\ph 
f_\ph(\tau)\,\vth(\tau,z,h;\ph)\tag1.42$$
defines an $SL_2(\Z)$ invariant function on $\frak H$. 
It can be very rapidly increasing on the standard fundamental domain for $\Gamma=SL_2(\Z)$. The regularization 
used to define the integral (1.39) will be reviewed in detail below. 

A basic result of Borcherds, \cite{\borch}, expressed in our present notation, is the following:
\proclaim{Theorem 1.2}{\rm(Theorem~13.3 of \cite{\borch})} Suppose that $F$ (and hence $f$) takes values in 
$S_M^K$ and that the Fourier coefficients $c_\ph(m)$ 
for $m\le 0$ are integers. Then the 
regularized integral
$$\P(z,h;F) = \II_{\Gamma\back \frak H}
(\!(\ f(\tau), \vth(\tau,z)\ )\!) \,v^{-2}\,du\,dv$$
can be written in the form
$$\P(z,h;F) = -2\log|\Psi(z,h;F)|^2 - c_0(0)\big(\ 2\log|y| + \log(2\pi) + \Gamma'(1)\ \big)$$
for a meromorphic modular form $\Psi(F)$ on $D\times H(\A_f)$ of weight $k=\frac 12c_0(0)$.
\endproclaim

More precisely, suppose that $c_0(0)$ is even, so that $k=\frac12 c_0(0)\in \Z$. 
Then, there is a unitary character $\xi$ of $H(\Q)$ such that, for all $\gamma\in
H(\Q)$, 
$$\Psi(\gamma z,\gamma h;F) = \xi(\gamma)\, j(\gamma,z)^k\,\Psi(z,h;F).\tag1.43$$ 
Moreover, as a function of $h\in H(\A_f)$, $\Psi(F)$ is right $K$--invariant for any compact open subgroup 
$K\subset H(\A_f)$ for which the values of $F$ lie in $S_M\subset S(V(\A_f))^K$, and hence, 
$\Psi(F)$ defines a meromorphic section of the bundle $\Cal L^{\tt k}\tt\Cal V_\xi$, where $\Cal V_\xi$ 
is the flat bundle defined by $\xi$. Since our calculations only involve $\log||\Psi(F)||^2$, the character $\xi$, 
which, in fact, has finite order \cite{\borchdukeII},  
will play no role in the present paper. If 
the coefficient $c_0(0)$ is odd, $\Psi(F)^2= \Psi(2F)$ is an automorphic form of weight $2k$. 
Note that, in any case, it is the quantity $2\log|\Psi(z,h;F)|^2$ which 
occurs in $\P(z,h;F)$, so that the parity of $c_0(0)$ will not matter. 

Borcherds also determines the divisor of $\Psi(F)$. To describe this in our setup, we first 
recall the definition of the special cycles in $X_K$, from \cite{\duke}. For $x\in V(\Q)$ with 
$Q(x)>0$, let $V_x = x^\perp$, and 
$$D_x=\{\ z\in D\mid x\perp z\ \}.\tag1.44$$
Let $H_x$ be the stabilizer of $x$ in $H$, and note that $H_x\simeq \GSpin(V_x)$. For 
$h\in H(\A_f)$, there is a natural map
$$\eqalign{
H_x(\Q)\back D_x\times H_x(\A_f)/(H(\A_f)\cap h K h^{-1}) &\lra H(\Q)\back D\times H(\A_f)/K = X_K\cr
\nass
(z,h_1)&\mapsto (z,h_1h)\cr}\tag1.45
$$
which defines a divisor $Z(x,h,K)$ on $X_K$. This divisor is rational over $\Q$. 
For a Schwartz function $\ph\in S(V(\A_f))^K$, and a positive rational number $m\in \Q_{>0}$, we define
a weighted linear combination $Z(m,\ph,K)$ of these divisors as follows. 
Let
$$\O_m=\{x\in V\mid Q(x)=m\}\tag1.46$$ 
be the quadric determined by $m$, and fix $x_0\in \O_m(\Q)$, assuming that $\O_m(\Q)\ne\phi$. 
Then $\O_m(\A_f)$ is a closed subset of $V(\A_f)$, and we can write
$$\text{\rm supp}(\ph)\cap \O_m(\A_f) = \coprod_r K\cdot \xi_r^{-1} x_0\tag1.47$$
for some finite set of $\xi_r$'s in $H(\A_f)$. 
Define
$$Z(m,\ph,K):= \sum_r\ph(\xi_r^{-1}x_0)\,Z(x_0,\xi_r,K).\tag1.48$$
If $\O_m(\Q)$ is empty, then $Z(m,\ph;K)=0$. 
These cycles, which are defined for arbitrary codimension in \cite{\duke}, include the Heegner points, 
Hirzebruch--Zagier curves, and Humbert surfaces as particular cases.  Various nice properties 
of the weighted cycles are described in \cite{\duke}. For example, if $K'\subset K$ and if $\pr:X_{K'}\rightarrow X_K$
is the associated covering, then
$$\pr^*Z(m,\ph,K) = Z(m,\ph,K'),\tag1.49$$
so that the special cycles are defined on the full Shimura variety associated to 
$(H,D)$, \cite{\milne}. Because of this relation, we will frequently omit $K$ and write
simply $Z(m,\ph)$ in place of $Z(m,\ph,K)$. 
Also, if $h\in H(\A_f)$, then right multiplication by $h^{-1}$ 
defines a natural morphism, rational over $\Q$,
$$r(h):X_K\lra X_{hKh^{-1}},\tag1.50$$
and
$$r(h)_*Z(m,\ph,K) = Z(m,\o(h)\ph,hKh^{-1}),\tag1.51$$
where $\o(h)\ph(x)=\ph(h^{-1}x)$. This relation describes the compatibility of 
the special cycles with the Hecke operators. Finally, by Proposition~5.4 of \cite{\duke},
we can give an explicit description of these cycles with respect to the decomposition (1.3) 
of the space $X_K$ as a disjoint union of arithmetic quotients of $D^+$:
$$Z(m,\ph;K)   = \sum_j \sum_{\matrix \scr x\in \O_m(\Q) \\\scr \mod \Gamma_j\endmatrix} \ph(h_j^{-1}x)\,\pr_j(D_x),\tag1.52$$
where $\pr_j:D^+\rightarrow \Gamma_j\back D^+$ is the natural projection. Note that it follows 
from this formula that,
$$Z(m,\ph;K) = Z(m,\ph^\vee;K)\tag1.53$$
where $\ph^\vee(x) = \ph(-x)$.

\proclaim{Theorem 1.3} {\rm (Theorem~13.3 of \cite{\borch})}
For $f$ with Fourier expansion given by (1.40) and (1.41), 
$$\div(\Psi(F)^2) =  \sum_\ph\sum_{m>0} c_\ph(-m)\,Z(m,\ph,K).$$
Here $\ph$ runs over the coset basis for $S_M$. 
\endproclaim
\demo{Proof} The singular part of the function $\P(F)$ is computed in \cite{\borch}. In our present notation, 
it is given as follows. Let $\b_1(t) = -\Ei(-t)$ be the exponential integral, \cite{\lebedev}, and 
recall that $\b_1(t) = -\log(t) + O(1)$ as $t\rightarrow 0$. The singular part of $\P(F)$ is given by
$$S(z,h;F) := \sum_{m>0} \sum_\ph c_\ph(-m)
\sum_{\scr x\in V(\Q)\atop\scr Q(x)=m} \ph(h^{-1}x)\,\b_1(-2\pi R(x,z)).\tag1.54$$
By (1.16), this has the same singularity as the logarithm of the absolute value of the 
infinite product
$$\prod_\ph\prod_{m>0} \prod_{\scr x\in V(\Q)\atop\scr Q(x)=m} (x,w(z))^{c_\ph(-Q(x))\ph(h^{-1}x)}.\tag1.55$$
Alternatively, in the neighborhood of any point $z_0\in \Bbb D$,  consider the divisor 
given by the finite product
$$\prod_\ph\prod_{m>0} \prod_{\matrix \scr x\in V(\Q) \\ 
\scr Q(x)=m\\ (x,w(z_0))=0\endmatrix } (x,w(z))^{c_\ph(-Q(x))\ph(h^{-1}x)}.\tag1.56$$
Note that $D_x$ is precisely the 
divisor defined by the equation $(x,w(z))=0$ on $\Bbb D$. 
Thus, since the singularities of $\P(F)$ on $\Bbb D$ are the same as those of $-2\log|\Psi(F)|^2$, we obtain
the claimed  result.
\qed\enddemo

\subheading{\Sec2. Computation of a regularized integral}

Setting $k=\frac12 c_0(0)$, recall that the Petersson norm of the section defined by $\Psi(F)$ is
$$||\Psi(z,h;F)||^2 = |\Psi(z,h;F)|^2 |y|^{2k}.\tag2.1$$
We view the function $||\Psi(z,h_j;F)||^2$ as a function on the component $\Gamma_j\back D^+$ 
of $X_K$ and will write $||\Psi(z;F)||^2$ for the resulting function on the (possibly disconnected) complex
manifold $X_K$. 

The basic problem is to compute the following integral:
$$\eqalign{
\kappa(\Psi(F)) &:= -\frac1{\vol(X_K)}\,\int_{X_K}\log||\Psi(z;F)||^2\, d\mu(z)\cr
\nass
\nass
{}&=-\frac1{\vol(X_K)}\,\int_{X_K}\log\big(\big|\Psi(z;F)\,\big|^2\, 
|y|^{2k}\big)\, d\mu(z)\cr
\nass
\nass
{}&= \frac12\, \frac1{\vol(X_K)}\,\int_{X_K}\P(z;F)\, 
d\mu(z) +k\, \big(\log(2\pi) +\Gamma'(1)\big)\cr 
\nass
\nass
{}&= \frac12\,\frac1{\vol(X_K)}\, \int_{X_K}
\bigg(\II_{\Gamma\back \frak H}
(\!(\ f(\tau), \vth(\tau,z)\ )\!)  \,v^{-2}\,du\,dv\,\bigg)\,  d\mu(z) +k\, C_0\cr
\nass
\nass
{}&= \frac12\,\II_{\Gamma\back \frak H} 
\sum_\ph f_\ph(\tau)\,I(\tau;\ph)\,v^{-2}\,du\,dv+k\,C_0,\cr}\tag2.2
$$
where $C_0=\log(2\pi) +\Gamma'(1)$, $\vol(X_K) = \vol(X_K, d\mu(z)\,dh)$ and 
$$I(\tau;\ph) = \frac1{\vol(X_K)}\,\int_{X_K}\vth(\tau,z;\ph)\,d\mu(z).\tag2.3$$
In fact, the last interchange of order of integration (where 
one of the integrals regularized!) will be 
justified in the next section, {\it provided} the theta integral (2.3) converges.
We will discuss this point in a moment. Here $d\mu(z)$ is a $H(\R)$--invariant top 
degree form on $D$; the quantity $\kappa(\Psi(F))$ is independent of the normalization 
of this form.

We want to relate the integral $I(\tau;\ph)$, over the complex manifold $X_K$, to 
the usual theta integral, over the 
adelic coset space $SO(V)(\Q)\back SO(V)(\A)$, appearing in the Siegel--Weil 
formula. This is done in detail in section 4, below, cf. Theorem~4.1. Note that there is 
an exact sequence
$$1\lra Z\lra H\lra SO(V)\lra1$$
where $H=\text{\rm GSpin}(V)$, as before. 
For simplicity, we assume that the compact open subgroup $K\subset H(\A_f)$ 
satisfies the condition: 
$$Z_K := Z(\A_f)\cap K \simeq \hat\Z^\times\tag2.4$$ 
under the natural identification $Z(\A_f)\simeq \A_f^\times$. A slight variant of the 
proof of Proposition~4.17 yields:
\proclaim{Lemma 2.1} Let $\ph_\infty$ be the Gaussian, as in (1.18) above. Then, for $\ph\in S(V(\A_f))^K$,
$$\align
I(g';\ph_\infty\tt\ph) &:= \int_{O(V)(\Q)\back O(V)(\A)} \theta(g',h;\ph_\infty\tt\ph)\, dh\\
\nass
\nass
{}&=\frac1{\vol(X_K)} \int_{X_K} \theta(g',z;\ph)\,d\mu(z).
\endalign
$$
\endproclaim
Note that both sides are independent of the choice of $d\mu(z)$. 
\proclaim{Corollary 2.2} 
$$I(\tau;\ph) = v^{-\ell/2}\,I(g'_\tau;\ph_\infty\tt\ph).$$
\endproclaim

\proclaim{Corollary 2.3} 
Assume that $F$ is valued in $S(V(\A_f))^{K}$. 
Then
$$\align
\kappa(\Psi(F)) &= \frac12\,\II_{\Gamma\back \frak H} 
\sum_\ph f_\ph(\tau)\,I(\tau;\ph)\,v^{-2}\,du\,dv+k\,C_0\\
\nass
\nass
{}&=\frac12\,\II_{\Gamma\back \frak H} 
\sum_\ph f_\ph(\tau)\,v^{-\ell/2}\,I(g'_\tau;\ph_\infty\tt\ph)\,v^{-2}\,du\,dv+k\,C_0,
\endalign
$$
with
$C_0 = \log(2\pi)+\Gamma'(1)$.
\endproclaim 

{\bf Remark 2.4.  Exceptional cases:} 
By Weil's criterion, \cite{\weilII}, p.75, Proposition 8, the theta integral $I(g'_\tau;\ph_\infty\tt \ph)$ is absolutely convergent 
whenever  $n-r>0$, where $r=0$, $1$, or $2$ is the Witt index 
of $V(\Q)$, i.e., the dimension of a maximal isotropic subspace of $V(\Q)$. Note that 
$r=0$ is only possible when $n\le 2$. The only {\it exceptional cases} will thus be 
$n=1$ with $V$ isotropic ($r=1$) and $n=2$ with $V$ split ($r=2$).  We will exclude these 
cases for now -- although they can be handled by the regularization process used in 
\cite{\krannals}.

We consider the regularized integral in the expression for $\kappa(\Psi(F))$ in Corollary 2.3.  Note that 
$\kappa(\Psi(F))$ is independent of the choice of the lattice $M$ and of $K$. 

Recall that, for a $\Gamma=\text{\rm PSL}_2(\Z)$ invariant function $\phi$ on $\frak H$, the regularized integral
$$\II_{\Gamma\back \frak H} \phi(\tau)\, d\mu(\tau),\tag2.5$$
used by Borcherds, 
is defined by taking the constant term in the Laurent expansion at $\s=0$ of the 
function defined, for $\Re(\s)$ sufficiently large, by
$$\lim_{T\rightarrow\infty} \int_{\Cal F_T} \phi(\tau)\,v^{-\s-2}\, du\, dv.\tag2.6$$
Here $\Cal F$ is the standard fundamental domain for the action of $\Gamma$ on 
$\frak H$, and $\Cal F_T$ is the intersection of this with the region $\Im(\tau)\le T$. 
This procedure can be applied provided that (i) the limit as $T$ goes to infinity exists in a halfplane $\Re(\s)>\s_0$, 
and (ii) the resulting holomorphic function of $\s$ has a meromorphic analytic continuation 
to a neighborhood of the point $\s=0$. In short, 
$$\II_{\Gamma\back \H} \phi(\tau)\ d\mu(\tau) := \CT{\s=0}\bigg\{\lim\limits_{T\rightarrow\infty}
\int_{\Cal F_T} \phi(\tau)\,v^{-\s}\,d\mu(\tau)\bigg\},\tag2.7$$
where $\CT{\s=0}$ denotes the constant term of 
the Laurent expansion at the point $\s=0$.  

The following result will be proved in 
the next section.

\proclaim{Proposition 2.5} 
$$\multline
\CT{\s=0}\bigg\{\lim\limits_{T\rightarrow\infty} \int_{\Cal F_T}\sum_\ph f_\ph(\tau)\,I(\tau;\ph)\, v^{-\s-2}\,du\,dv\bigg\}\\
\nass
{}= \lim\limits_{T\rightarrow \infty}\bigg[\  \int_{\Cal F_T}\sum_\ph f_\ph(\tau)\,I(\tau;\ph)\, v^{-2}\,du\,dv
-c_0(0)\,\log(T)\ \bigg].
\endmultline
$$
\endproclaim

Thus we need to evaluate the basic integral
$$ \int_{\Cal F_T}\sum_\ph f_\ph(\tau)\,I(\tau;\ph)\, v^{-2}\,du\,dv\tag2.8$$
where $f_\ph(\tau)$ is holomorphic on $\Cal F$ and where $\s$ has been set equal to zero. 

Following the suggestion of section 9 of \cite{\borch}, we would {\it like} to define an automorphic function 
$J(\tau;\ph)$ on $\frak H$ for which 
$$\frac{\partial}{\partial \bar\tau}\big\{J(\tau;\ph)\big\} = I(\tau;\ph)\,v^{-2}.\tag2.9$$
Then, by a simple Stokes' Theorem argument, we would have
$$\eqalign{
&\int_{\Cal F_{T}} \sum_\ph f_\ph(\tau)\,I(\tau;\ph)\, v^{-2}\,du\wedge dv\cr
\nass
{}&=\frac1{2i}\int_{\Cal F_{T}} d\bigg(\sum_\ph f_\ph(\tau)\, J(\tau;\ph)\,d\tau\bigg)\cr
\nass
{}&=\frac1{2i}\int_{\d \Cal F_{T}} \sum_\ph f_\ph(\tau)\, J(\tau;\ph)\,d\tau\cr 
\nass
{}&=\frac1{2i}\int_{1/2+iT}^{-1/2+iT} \sum_\ph f_\ph(\tau)\, J(\tau;\ph)\,du\cr
\nass
{}&= -\frac1{2i}\ \text{constant term of}\ \bigg( \sum_\ph f_\ph(\tau)\, J(\tau;\ph)\bigg)\bigg|_{v=T}.\cr
}\tag2.10
$$
In the next to last step, we have used the invariance of $\sum_\ph f_\ph(\tau)\, J(\tau;\ph)\,d\tau$ 
under $\tau\mapsto \tau+1$ and under $\tau\mapsto -1/\tau$. 

To obtain a relation like (2.9), we apply the Maass operators and the Siegel--Weil formula. 
Let
$$X_{\pm} = \frac12 \pmatrix 1&\pm i\\\pm i&-1\endpmatrix \in \frak{sl}_2(\C).\tag2.11$$
Recall that, if $\phi:G'_\R\rightarrow \C$ is a smooth function with $\phi(g'k') = \chi_{\ell}(k')\phi(g')$, 
i.e., of weight $\ell$,  
and if $f(\tau) = v^{-\frac{\ell}2}\,\phi(g'_\tau)$ is the corresponding function on $\frak H$, 
then $X_\pm \phi$ has weight $\ell \pm 2$, and the corresponding function on $\frak H$ is 
$$ v^{-\frac12(\ell\pm2)}X_{\pm}\phi(g'_\tau) =\cases \big( 2 i \frac{\d f}{\d\tau} + \frac{\ell}{v}f\big)(\tau) 
&\text{ for $+$,}\\
\nass
\nass
-2i v^{2} \ \frac{\d f}{\d \bar\tau}(\tau)&\text{ for $-$.}
\endcases\tag2.12
$$

We now take advantage of the 
Siegel--Weil formula; the facts we need are reviewed in the first part of section 4.  
For $\ph\in S(V(\A_f))$, 
let $E(g',s,\P^r_\infty\tt\lambda(\ph))$
be the Eisenstein series 
of weight $r$ on $G'_\A$ associated to $\ph$.  If $\ph_\infty\in S(V(\R))$ is the Gaussian, then 
$$\lambda(\ph_\infty) =\P^\ell_\infty(s_0),\tag2.13$$
where $\ell =\frac{n}2-1$, as above. By the Siegel--Weil formula, Theorem~4.1, we have the following. 
\proclaim{Proposition 2.6} Exclude the exceptional cases of Remark~2.4 above, so that the 
theta integral is 
absolutely convergent. Then
$$I(\tau;\ph) = v^{-\ell/2}\,I(g'_\tau;\ph_\infty\tt\ph) = v^{-\ell/2}\,E(g'_\tau,s_0;\P^\ell_\infty\tt\lambda(\ph)),$$
where $s_0 = \frac{n}2=\ell+1$. 
\endproclaim

On the other hand, an easy computation in the induced representation $I_\R(s,\chi)$ of $G'_\R$ shows:
\proclaim{Lemma 2.7} Let $\P^r(s)\in I_\R(s,\chi)$ be the normalized eigenvector of weight $r$ for the action 
of $K_\R'$. Then
$$X_\pm \P^r(s) = \frac12 (s+1\pm r)\, \P^{r\pm2}(s).$$
\qed
\endproclaim
Therefore, we have the basic relation:
$$X_- E(g',s;\P^{\ell+2}\tt\lambda(\ph)) = \frac12(s-\ell-1)\, E(g',s;\P^{\ell}\tt\lambda(\ph)).\tag2.14$$
Pushing this down to $\frak H$ we obtain:
$$-2i v^2 \ \frac{\d }{\d \bar\tau}\bigg\{ v^{-\frac12(\ell+2)}\, E(g'_\tau,s;\P^{\ell+2}\tt\lambda(\ph))\bigg\} 
= \frac12(s-s_0)\,v^{-\frac12\ell} E(g'_\tau,s;\P^{\ell}\tt\lambda(\ph)).\tag2.15$$
For convenience, we now write
$$E(\tau,s;\ph,\ell)=v^{-\ell/2}\,E(g'_\tau,s_0;\P^\ell_\infty\tt\lambda(\ph)),\tag2.16$$
so that (2.15) becomes
$$-2i v^2 \ \frac{\d }{\d \bar\tau}\bigg\{  E(\tau,s;\ph,\ell+2)\bigg\} 
= \frac12(s-s_0)\, E(\tau,s;\ph,\ell).\tag2.17$$
Of course,the vanishing of the right hand side of (2.17) at $s=s_0= \frac{n}2$ just shows the holomorphy of
the special value
$$ E(\tau,s_0;\ph,\ell+2)
= \ph(0) + \vol(X)^{-1} \sum_{m>0} \deg(Z(m,\ph))\cdot q^m,\tag2.18$$
cf. Theorem~4.23. Here we have written $\deg(Z(m,\ph))$ 
in place of 
$$\deg_{\Cal L^\vee}(Z(m,\ph;K))$$
and $\vol(X)$ in place of $\vol(X_K,\O^n)$ to lighten the notation. 

{\bf Remark 2.8.} The vanishing of the right side of (2.17) depends on the fact that 
$E(\tau,s;\ph,\ell)$ has no pole at 
$s=s_0=\frac{n}2$. In the exceptional cases, $n=1$, $r=1$ and $n=2$, $r=2$ a pole can occur, and
its residue accounts for a non--holomorphic 
component occuring in (2.14), cf. \cite{\funke}.

We write:
$$E(\tau,s;\ph,\ell)\, v^{-2} 
= \frac{-4i}{s-s_0}\, \frac{\d }{\d \bar\tau}\bigg\{  E(\tau,s;\ph,\ell+2)\bigg\}.\tag2.19$$
Now, to evaluate (2.8), we use the Siegel--Weil formula (Proposition~2.6) and write
$$\int_{\Cal F_T}\sum_\ph f_\ph(\tau)\,I(\tau;\ph)\, v^{-2}\,du\,dv = 
\int_{\Cal F_{T}} \sum_\ph f_\ph(\tau)\, E(\tau,s;\ph,\ell)\, v^{-2}\,du\wedge dv\,\bigg|_{s=s_0}.$$
Then, for general $s$, we can use the relation (2.19) and the Stoke's Theorem argument (2.10) to obtain the following basic identity.
$$\eqalign{
I(s,T):=&\int_{\Cal F_{T}} \sum_\ph f_\ph(\tau)\,
E(\tau,s;\ph,\ell)\, v^{-2}\,du\wedge dv\cr
\nass
\nass
{}&=\frac1{2i}\int_{\Cal F_{T}} d\bigg(\sum_\ph f_\ph(\tau)\,
\frac{-4i}{s-s_0}   E(\tau,s;\ph,\ell+2)
 \,d\tau\bigg)\cr
\nass
\nass
{}&=\frac{-2}{s-s_0} \int_{\d \Cal F_{T}} 
\sum_\ph f_\ph(\tau)\, E(\tau,s;\ph,\ell+2)\,d\tau\cr
\nass
\nass
{}&=\frac{-2}{s-s_0} \int_{1/2+iT}^{-1/2+iT} 
\sum_\ph f_\ph(\tau)\, E(\tau,s;\ph,\ell+2)\,du\cr
\nass
\nass
{}&= \frac{2}{s-s_0} \cdot \text{constant term of}\ 
\bigg( \sum_\ph f_\ph(\tau)\, E(\tau,s;\ph,\ell+2)\bigg)\bigg|_{v=T}.\cr}
\tag2.20
$$
By Corollary~2.3, and Proposition~2.5,
$$\align
\kappa(\Psi(F)) &= \frac12\,\lim\limits_{T\rightarrow \infty}\bigg[\  \int_{\Cal F_T}\sum_\ph f_\ph(\tau)\,I(\tau;\ph)\, v^{-2}\,du\,dv
-c_0(0)\,\log(T)\ \bigg]+k\,C_0\\
\nass
\nass
{}&=\frac12\,\lim\limits_{T\rightarrow \infty}\bigg[\  I(s_0,T)
-c_0(0)\,\log(T)\ \bigg]+k\,C_0,
\endalign
$$
It will be convenient to introduce the following additional notation. 
Write
$$E(\tau,s;\ph,\ell+2) = \sum_m A_\ph(s,m,v)\, q^m,\tag2.22$$ 
where the Fourier coefficients have Laurent expansions 
$$A_\ph(s,m,v) = a_\ph(m) + b_\ph(m,v)(s-s_0) + O((s-s_0)^2).\tag2.23$$
where the $a_\ph(m)$'s are given by (2.16). 
With this notation, 
$$\eqalign{
I(s,T)&=\frac{2}{s-s_0}\cdot \text{constant term of}\ 
\bigg(\sum_\ph f_\ph(\tau)\, E(\tau,s;\ph,\ell+2)\bigg)\bigg|_{v=T}\cr
\nass
\nass
&{}\qquad=\frac{2}{s-s_0}\sum_\ph\sum_m c_\ph(-m)\, A_\ph(s,m,T).\cr}\tag2.24
$$
We consider the individual terms. 
For $m=0$, we have
$$\frac2{s-s_0}\sum_\ph c_\ph(0)\, \big(\, \ph(0) + b_\ph(0,T) (s-s_0)\,\big)  + O(s-s_0).\tag2.25$$
so that the contribution of such terms to the constant coefficient in the Laurent expansion at $s=s_0$ is 
$$2 \sum_\ph c_\ph(0)\, b_\ph(0,T).\tag2.26$$
We will return to the polar part occuring in (2.25) in a moment. 
Similarly, from the $m<0$ terms, we have the contribution
$$2\sum_\ph \sum_{m<0} c_\ph(-m)\, b_\ph(m,T).\tag2.27$$ 
Finally, for the finite sum of terms with $m>0$, we have, initially:
$$\eqalign{
\frac{1}{(s-s_0)}\frac{2}{\vol(X)}&\sum_\ph \sum_{m>0} c_\ph(-m)\,\deg(Z(m,\ph))\cr
\nass
 {}&+ 2 \sum_\ph\sum_{m>0} c_\ph(-m)\,b_\ph(m,T) + O(s-s_0).\cr}\tag2.28$$ 
Since our whole integral $I(s,T)$ does not have a pole at $s=s_0$, the polar part here must cancel the 
one which occurred earlier, i.e., we must have
$$2\sum_\ph c_\ph(0)\, \ph(0) +\frac{2}{\vol(X)} \sum_\ph \sum_{m>0} c_\ph(-m)\,\deg(Z(m,\ph)) = 0.\tag2.29$$
Since
$$\div(\Psi(F)^2) = \sum_{m>0}c_\ph(-m)\,Z(m,\ph),\tag2.30$$
this amounts to
$$\deg(\div(\Psi(F)^2) = \sum_{m>0}c_\ph(-m)\,\deg(Z(m,\ph)) = -\vol(X)\,c_0(0).\tag2.31$$
Recall that we are using the coset basis for $S_M$, so that $\ph_0(0)=1$ and $\ph(0)=0$ for 
$\ph\ne \ph_0$. Also note that, since $\O$ is the negative of a K\"ahler form, $\vol(X)$ 
and $\deg(Z(m,\ph))$ will have opposite signs (for a coset function $\ph$), cf. (4.49).

{\bf Example 2.9.} Suppose that $n=1$ and $r=0$, i.e., 
$V$ is anisotropic over $\Q$ of dimension $3$ and $X_K$ is a (disjoint union of) 
projective curves. Suppose that the image of $K$ in $SO(V)(\A_f)$ is neat, so that all of the $\Gamma_j$'s 
act without fixed points on $D^+\simeq \frak H$.  Then, 
since $\O = -\frac{1}{2\pi}\,y^{-2}\,dx\wedge dy,$ 
$\vol(X) = 2-2g$, where $g$ is the genus of $X_K$, and hence we have
$$\deg(\div(\Psi(F)^2) = 2(g-1) \,c_0(0),\tag2.32$$
as expected. Here one must keep in mind the fact that $\Psi(F)^2$ has 
`classical weight' $2\,c_0(0)$. 

Collecting the contributions of (2.26), (2.27), and (2.28), we obtain
\proclaim{Proposition 2.10}
$$\eqalign{
I(s_0,T) &= \int_{\Cal F_{T}}\sum_\ph f_\ph(\tau)\,I(\tau,\ph)\, v^{-2}\, du\, dv\cr
\nass
{}&=  2 \sum_\ph\sum_m c_\ph(-m)\,b_\ph(m,T).\cr}
$$
\endproclaim
The following result will be proved in the next section. 
\proclaim{Proposition 2.11}  (i) For $m<0$, $b_\ph(m,T)$ decays exponentially as $T\rightarrow\infty$. 
\hfill\break
(ii) 
$$\lim\limits_{T\rightarrow\infty}\bigg(\ 2 \sum_\ph\sum_{m<0} c_\ph(-m)\,b_\ph(m,T)\ \bigg)
 =0.$$
(iii) For $m=0$, 
$$\lim\limits_{T\rightarrow\infty}\bigg(\, b_0(0,T) - \frac12\,\log(T) \,\bigg) =0, $$
and, for $\ph\ne\ph_0$, 
$$\lim\limits_{T\rightarrow\infty} b_\ph(0,T) = 0.$$
\endproclaim

Thus, we obtain an explicit expression for the quantity $\kappa(\Psi(F))$. The following result 
summarizes the relations between the geometry of the Borcherds form $\Psi(F)$ and the 
family of Eisenstein series $E(\tau,s,;\ph,\ell+2)$. 

\proclaim{Main Theorem 2.12} For $\ph\in S(V(\A_f))$, let
$$E(\tau,s;\ph,\ell+2) = \sum_m A_\ph(s,m,v)\,q^m,$$
with
$$A_\ph(s,m,v) = a_\ph(m) + b_\ph(m,v)\,(s-s_0) + O((s-s_0)^2)$$
be the Laurent expansion at the point $s_0= \frac{n}2 = \ell+1$ of the 
associated Eisenstein series of weight $\frac{n}2+1 = \ell+2$. Let $K$ be a  
compact open subgroup $K\subset H(\A_f)$ satisfying the condition (2.4), and
let $X=X_K$.
Exclude the cases $\dim V=3$, of Witt index $1$ and $\dim V=4$, of Witt index $2$.\hfill\break
(i) Suppose that $\ph\in S(V(\A_f))^{K}$. Then, 
$$E(\tau,s_0;\ph,\ell+2) = \ph(0) + \vol(X)^{-1}\sum_{m>0}\deg_{\Cal L^\vee}(Z(m,\ph))\,q^m.$$
(ii) For any $\ph\in S(V(\A_f))$, let
$$\kappa_\ph(m) :=\cases \lim\limits_{T\rightarrow\infty} b_\ph(m,T) &\text{ if $m>0$, and}\\
\nass
\frac12 C_0\,\ph(0)&\text{ if $m=0$,}
\endcases
$$
where $C_0=\log(2\pi)+\Gamma'(1)$. Suppose that $f:\frak H \rightarrow S(V(\A_f))^{K}$
is a modular form of weight $1-\frac{n}2 = -\ell$ for $\SL_2(\Z)$, 
with Fourier expansion
$$f(\tau) = \sum_{\ph}\sum_m c_\ph(m)\,q^m\,\ph$$
where $\ph$ runs over the coset basis with respect to some lattice $M$ and where
$c_\ph(m)\in\Z$ for $m\le0$. Let $\Psi(f)$ be the associated Borcherds form of weight 
$c_0(0)/2$. 
Then
$$\div(\Psi(f)^2) = \sum_{\ph}\sum_{m>0}c_\ph(-m)\,Z(m,\ph),$$
and
$$-\vol(X)\,c_0(0) = \sum_{\ph}\sum_{m>0}c_\ph(-m)\,\deg_{\Cal L^\vee}(Z(m,\ph)).$$
Moreover
$$\eqalign{
\kappa(\Psi(f)):&= -\frac{1}{\vol(X)} \int_{X_K}\log||\Psi(z;f)||^2\,d\mu(z)\cr
\nass
\nass
{}&= \sum_\ph\sum_{m\ge0} c_\ph(-m)\,\kappa_\ph(m).\cr}
$$
Here $\vol(X) = \vol(X_K,\O^n)$ and $\deg_{\Cal L^\vee}(Z(m,\ph)) = \int_{Z(m,\ph;K)} \O^{n-1}$
are computed with respect to the invariant $(1,1)$--form $\O =dd^c\log(\rho)$, where $\rho = \rho(z) = -\frac12(w(z),w(\bar z))$, 
cf. Proposition~4.10. 
\endproclaim

{\bf Remark 2.13:} The quantity $\kappa(\Psi(f))$ is completely determined by the 
collection of integers $\{\,c_\ph(-m)\,\}$ for $m\ge0$. The universal quantities $\kappa_\ph(m)$ 
are independent of $\Psi(f)$. 
They can be computed explicitly, cf. section 5 for an example and \cite{\ky} for a more systematic
discussion.

\subheading{\Sec3. Convergence estimates}

In this section we prove the crucial fact that the integration over $X_K$ 
can be interchanged with the Borcherds' regularization. 
\proclaim{Theorem 3.1} Suppose that the integral of the theta function converges, 
i.e., suppose that $V$ is not a ternary isotropic space of signature $(1,2)$ 
or a quaternary space of signature $(2,2)$ and $\Q$--rank 2, the 
exceptional cases of Remark~2.4 above. Then
$$\int_{X_K} \II_{\Gamma\back \H} ((\ F(\tau),\vth(\tau,z)\ ))\ d\mu(\tau)\,d\mu(z),
=\II_{\Gamma\back \H}  ((\ F(\tau),\int_{X_K}\vth(\tau,z) \,d\mu(z)\ )) \, d\mu(\tau),$$
where $\II$ denotes the regularized integral. 
\endproclaim

Writing
$\Cal F_T = \Cal F_1 \cup \Cal B_T$, where $\Cal B_T= \Cal F_T-\Cal F_1$, 
we consider the first expression: 
$$\align
&\int_{X_K} \II_{\Gamma\back \H}  ((\ F(\tau),\vth(\tau,z)\ ))\ d\mu(\tau)\,d\mu(z)\\
\nass
{}&= \int_{X_K} \CT{\s=0}\bigg\{\,\lim\limits_{T\to\infty}
\int_{\Cal F_T} ((\ F(\tau),\vth(\tau,z)\ ))\, v^{-\s}\,d\mu(\tau)\,\bigg\}
\,d\mu(z)\\
\nass
\nass
{}&=\int_{X_K} \CT{\s=0}\bigg\{\,\lim\limits_{T\to\infty}
\int_{\Cal B_T} ((\ F(\tau),\vth(\tau,z)\ ))\, v^{-\s}\,d\mu(\tau)\,\bigg\}
\,d\mu(z)\tag3.1\\
\nass
{}&\qquad\qquad + \int_{X_K}\int_{\Cal F_1} ((\ F(\tau),\vth(\tau,z)\ ))\, d\mu(\tau)\,d\mu(z)\\
\nass
\nass
{}&=\int_{X_K} \CT{\s=0}\bigg\{\,\lim\limits_{T\to\infty}\int_{1}^T
C(v,z)\, v^{-\s-1}\,dv\,\bigg\}
\,d\mu(z)\\
\nass
{}&\qquad\qquad + \int_{\Cal F_1} \int_{X_K}((\ F(\tau),\vth(\tau,z)\ ))\,d\mu(z)\, d\mu(\tau),
\endalign
$$
where 
$$\eqalign{
C(v,z) =C(v,z,h)&:= v^{-1}\int_{-\frac12}^{\frac12} ((\ F(\tau),\vth(\tau,z,h)\, ))\,du\cr
\nass
\nass
{}&=  \sum_\ph\sum_{m\in \Q} c_\ph(-m) \sum_{\scr x\atop \scr Q(x)=m}\ph(h^{-1}x)\, e^{-2\pi v R(x,z)}\cr
}\tag3.2
$$
is the constant term of  $v^{-1}((\ F(\tau),\vth(\tau,z)\, ))$.  
Here, in the term arising from integration over $\Cal F_1$, we have used the integrability of $\vth(\tau,z)$ over $X_K$. 
It now suffices to show that 
the term 
$$A:=\int_{X_K} \CT{\s=0}\bigg\{\,\lim\limits_{T\to\infty}\int_{1}^T
C(v,z)\, v^{-\s-1}\,dv\,\bigg\}
\,d\mu(z)\tag3.3$$
in the last expression can be rewritten as
$$B:=\CT{\s=0}\bigg\{\,\lim\limits_{T\to\infty}\int_{X_K} \int_{1}^T
C(v,z)\, v^{-\s-1}\,dv
\,d\mu(z)\,\bigg\}.\tag3.4
$$
To see this, observe that the integral in $B$ is then equal to 
$$\align
{}&\int_{X_K} \int_{1}^T C(v,z)\, v^{-\s-1}\,dv\,d\mu(z)\\
\nass
{}&=\int_{X_K}\int_{\Cal B_T}  ((\ F(\tau),\vth(\tau,z)\ ))\, v^{-\s}\,d\mu(\tau)\,d\mu(z)\tag3.5\\
\nass
{}&=\int_{\Cal B_T} \int_{X_K} ((\ F(\tau),\vth(\tau,z)\ ))\,d\mu(z)\, v^{-\s}\,d\mu(\tau),\\
\endalign
$$
again using the integrability of $\vth(\tau,z)$.  
Substituting the 
resulting expression for $B$ in place of $A$ in the last expression of (3.1), we obtain
$$
\align
&\CT{\s=0}\bigg\{\,\lim\limits_{T\to\infty}\int_{\Cal B_T} \int_{X_K} ((\ F(\tau),\vth(\tau,z)\ ))\,d\mu(z)\, v^{-\s}\,d\mu(\tau)\,\bigg\}\\
\nass
{}&\qquad\qquad + \int_{\Cal F_1} \int_{X_K}((\ F(\tau),\vth(\tau,z)\ ))\,d\mu(z)\, d\mu(\tau)\tag3.6\\
\nass
\nass
{}&=\CT{\s=0}\bigg\{\,\lim\limits_{T\to\infty}\int_{\Cal F_T} \int_{X_K} ((\ F(\tau),\vth(\tau,z)\ ))\,d\mu(z)\, v^{-\s}\,d\mu(\tau)\,\bigg\}\\
\nass
\nass
{}&=\II_{\Gamma\back\H} \int_{X_K} ((\ F(\tau),\vth(\tau,z)\ ))\,d\mu(z)\,
\endalign
$$
as required. 

To show the equality of $A$ and $B$, we break the function $C(v,z)$ into pieces. 
$$\align
C_+(v,z)&:= \sum_\ph\sum_{m>0} c_\ph(-m) \sum_{\scr x\atop \scr Q(x)=m}\ph(x)\, e^{-2\pi v R(x,z)}\\
\nass
C_0(v,z)&:= \sum_\ph c_\ph(0) \sum_{\scr x\atop \scr Q(x)=0, x\ne0}\ph(x)\, e^{-2\pi v R(x,z)}\tag3.7\\
\nass
C_{00}(v,z)&:= \sum_\ph c_\ph(0) \ph(0)  = c_0(0)\qquad\text{(for the coset basis)}\\
\nass
C_-(v,z)&:= \sum_\ph\sum_{m<0} c_\ph(-m) \sum_{\scr x\atop \scr Q(x)=m}\ph(x)\, e^{-2\pi v R(x,z)}.
\endalign
$$  
We will write $A_+$, $A_0$, $A_{00}$, and $A_-$
(resp. $B_+$, etc.) for the corresponding contributions to $A$ (resp. $B$). 

For the $C_{00}$ term, we have
$$
\int_1^T v^{-\s-1}\,dv = \frac1{\s}(1-T^{-\s}),\tag3.8$$
so that
$$\CT{\s=0}\bigg\{\,\lim\limits_{T\to\infty} \int_1^T C_{00}(v,z)\,dv \,\bigg\} =0.\tag3.9$$
This gives $A_{00}=B_{00}=0$. 

Next consider the quatities $A_+$ and $B_+$ arising from $C_+(v,z)$.
Note that the sum on $m>0$ in $C_+(v,z)$ is finite, since there 
are only finitely many nonvanishing negative Fourier coefficients $c_\ph(-m)$. 
For a given coset representative $h=h_j$, we write $\Gamma = \Gamma_j = H(\Q)\cap hKh^{-1}$, so that 
$\Gamma\back D^+$ is the associated component of $X_K$. For a fixed $m>0$ and $\ph$ 
and on the chosen component of $X_K$, 
the sum in $C_+(v,z)$  involves 
$$\{ x\in V(\Q)\mid Q(x) = m, \ \ph(h^{-1}x) \ne 0\, \}.\tag3.11$$
This set consists of a finite number of $\Gamma$ orbits. The contribution to $A$ 
of a single such orbit is
$$c_\ph(-m)\,\ph(h^{-1}x)\int_{\Gamma\back D^+} \CT{\s=0}\bigg\{\,\lim\limits_{T\to\infty}\int_{1}^T
\sum_{\gamma\in \Gamma_x\back \Gamma}\, 
e^{-2\pi v R(x,\gamma z)}\, v^{-\s-1}\,dv\,\bigg\}
\,d\mu(z).\tag3.12$$
To prove the finiteness of this expression, it will suffice to prove the finiteness of 
$$
\int_{\Gamma\back D^+} \lim\limits_{T\to\infty}\int_{1}^T
\sum_{\gamma\in \Gamma_x\back \Gamma}\, 
e^{-2\pi v R(x,\gamma z)}\, v^{-\s-1}\,dv\,
\,d\mu(z),\tag3.13$$
for $\s=\s_0$ for some real $\s_0<0$.  Indeed, such finiteness implies that (3.13) 
defines a holomorphic function of $\s$ in the half plane $\Re(\s)>\s_0$. If $z$ lies in the set 
$$ D - \bigcup_{\gamma\in \Gamma_x\back \Gamma} \gamma^{-1}D_x,\tag3.14$$
then none of the $R(x,\gamma z)$'s vanish and the limit on $T$ inside the 
integral is finite. Note that the excluded set of $z$'s has measure zero. The following result 
will be proved at the end of this section. 
\proclaim{Proposition 3.2} Let
$$\beta_{\s+1}(t) = \int_1^\infty e^{-t v}\, v^{-\s-1}\,dv.$$
Then, if $Q(x)>0$, the integral 
$$\align
&\int_{\Gamma\back D^+} \lim\limits_{T\to\infty}\int_{1}^T
\sum_{\gamma\in \Gamma_x\back \Gamma}\, 
e^{-2\pi v R(x,\gamma z)}\, v^{-\s-1}\,dv\,
\,d\mu(z)\\
\nass
\nass
{}&=\int_{\Gamma\back D^+} 
\sum_{\gamma\in \Gamma_x\back \Gamma}\, 
\beta_{\s+1}(2\pi R(x,\gamma z)) 
\,d\mu(z)\\
\nass
\nass
{}& = 
\int_{\Gamma_x\back D^+} \beta_{\s+1}(2\pi R(x,z)) 
\,d\mu(z)
\endalign
$$
is holomorphic in the halfplane $\Re(\s)> -1$.
\endproclaim
Recall that $\b_1(t) = O(-\log(t))$ as $t\rightarrow 0$ and 
$\b_1(t) = O(e^{-t})$ as $t\rightarrow\infty$. 
Thus, when $\s=0$, the integrand 
$\beta_1(2\pi R(x,z))$ has a logarithmic singularity on the `waist' 
$\Gamma_x\back D^+_x$ of the tube $\Gamma_x\back D^+$. Also note that 
this `waist' can be noncompact. 

\proclaim{Corollary 3.3} $A_+ = B_+$. 
\endproclaim

Next we consider the terms $A_0$ and $B_0$ associated to the nonzero null vectors. 
Again for a given $h$ and $\ph$, the associated terms in $C_0(v,z)$ will be 
$$c_\ph(0) \sum_{\scr x \ne0\atop \scr Q(x)=0}\ph(h^{-1}x)\, e^{-2\pi v R(x,z)}.\tag3.15$$
There are a finite number of $\Gamma$ orbits in the space of null lines in $V(\Q)$. 
For a given null line $\ell\subset V$, we have the contribution to $A_0$:
$$c_\ph(0)\int_{\Gamma\back D^+} \CT{\s=0}\bigg\{\,\lim\limits_{T\to\infty}\int_{1}^T
\sum_{\gamma\in \Gamma_\ell\back \Gamma}\, 
\sum_{x\in \ell(\Q), \ x\ne0}
\ph(h^{-1}x)\,e^{-2\pi v R(x,\gamma z)}\, v^{-\s-1}\,dv\,\bigg\}
\,d\mu(z).\tag3.16$$
Again, the following result, to be proved below, will suffice. 
\proclaim{Proposition 3.4} Suppose that $n>1$. Then the integral
$$
\int_{\Gamma_\ell\back D^+} 
\sum_{x\in \ell(\Q), \ x\ne0}\ph(h^{-1}x)\,
\b_{\s+1}(2\pi v R(x,z))\,
\,d\mu(z)$$
is holomorphic in the halfplane $\Re(\s)>-\frac{n}2$. 
\endproclaim
\proclaim{Corollary 3.5} $A_0=B_0$. 
\endproclaim

Finally, we turn to the terms where $m<0$. Note that the sum on $m$ in $C_-(v,z)$ 
is now infinite so that we will need information about the growth of the 
Fourier coefficients $c_\ph(-m)$. In fact, these can grow very fast!

As before, we fix $\ph$ and $h$, and, taking the limit with respect to $T$, we consider
$$\int_{\Gamma\back D^+}\sum_{m<0} c_\ph(-m) 
\sum_{\scr x\atop \scr Q(x)=m}\int_1^\infty \ph(h^{-1}x)\, e^{-2\pi v R(x,z)}\,v^{-\s-1}\,dv\,d\mu(z).\tag3.17$$
Here we can push the integral over $\Gamma\back D^+$ inside the sum on $m$, 
and again use the fact that, for each $m$,  
there are only a finite number of $\Gamma$ orbits in the set 
$$\{ x\in V(\Q)\mid Q(x)=m,\ \ph(h^{-1}x)\ne0\ \}.\tag3.18$$ 
Thus, it will suffice to show:
\proclaim{Proposition 3.6} The sum 
$$\sum_{m<0} c_\ph(-m) 
\sum_{\scr x\atop {\scr Q(x)=m \atop \scr \mod \Gamma}} \ph(h^{-1}x)\,
\int_{\Gamma_x\back D^+}\int_1^\infty  e^{-2\pi v R(x,z)}\,v^{-\s-1}\,dv\,d\mu(z)$$
defines an entire function of $\s$. 
\endproclaim

\proclaim{Corollary 3.7}
$A_-=B_-$. 
\endproclaim

\demo{\bf Proof of Proposition 3.2 } 

To show the finiteness of the integral 
$$\int_{\Gamma_x\back D^+} \beta_{\s+1}(2\pi R(x,z)) 
\,d\mu(z)\tag3.19$$
in the case $Q(x)>0$, we introduce coordinates. 
We choose a basis for $V(\R)$ so that the inner product has matrix $I_{n,2}$ 
and so that $x = 2\a v_1$ is a nonzero multiple of the first basis vector. 
Then $SO(V)(\R)^+\simeq SO^+(n,2)=G$ and the subgroup stabilizing $x$ is isomorphic to 
$SO^+(n-1,2)=G_x$. 
Let $z_0\in D^+$ be the oriented negative $2$--plane spanned by $v_{n+1}$ and $v_{n+2}$ and let $K=SO(n)\times SO(2)$ 
be its stabilizer in $SO^+(n,2)$. 
The plane spanned by $v_1$ and $v_{n+1}$, the first negative basis vector,
has signature $(1,1)$. The identity component of the special orthogonal group of this plane 
is a 1-parameter subgroup
$$A=\{a_t\mid t\in \R\},\tag3.20$$
where $a_t v_1 = \cosh(t)v_1+\sinh(t) v_{n+1}$. Let $A_+$ be the subset of $a_t$'s with $t\ge 0$. 
Then, from the general theory of 
semisimple symmetric spaces -- a convenient reference is \cite{\flenstedj} -- one has a double coset decomposition
$$G = G_xA_+K\tag3.21$$ 
and the integral formula
$$\int_{G} \phi(g) \, dg = \int_{G_x}\int_{A_+}\int_K \phi(g_x a_t k)\, 
|\sinh(t)|\cosh(t)^{n-1}\,dg_x\,dt\,dk.\tag3.22$$

For $z= g_xa_t\cdot z_0\in D^+$, we have
$$R(x,z) = 2m\sinh^2(t),\tag3.23$$
since $Q(x) = 2 \a^2 = m$. Then, our integral becomes (up to a positive constant 
depending on normalization of invariant measures)
$$\align
&\int_{\Gamma_x\back D^+}  \b_{\s+1}(2\pi R(x,z))\,d\mu(z)\tag3.24\\
\nass
{}&=C\, \vol(\Gamma_x\back G_x)\,\vol(K)
\int_0^\infty \b_{\s+1}(4\pi m\,\sinh^2(t))\,\sinh(t)\cosh(t)^{n-1}\,dt.
\endalign
$$ 

\proclaim{Lemma 3.8} (i) The function
$$\b_{\s+1}(t) = \int_1^\infty e^{-tu}\,u^{-\s-1}\, du$$
is $O(e^{-t})$ as $t\rightarrow\infty$. \hfill\break
(ii) If $\s<0$, then $\b_{\s+1}(t) = O(t^\s)$ as $t\rightarrow0$.
\hfill\break
(iii) If $\s=0$, then 
$$\b_1(t) = -\Ei(-t) = 
 -\log(t) + \gamma + \int_0^t \frac{e^u-1}{u}\,du
$$
is the exponential integral and this function has a logarithmic 
singularity, $-\log(t)$,  as $t\rightarrow 0$. 
\hfill\break
(iv) If $\s>0$, then $\b_{\s+1}(t) = O(1)$ as $t\rightarrow0$. \qed
\endproclaim

The integral (3.24) is finite for $\s>-1$, 
since, near the lower endpoint it looks like
$$\int_0 \sinh(t)^{2\s}\,\sinh(t)\cosh(t)^{n-1}\,dt = \int_0 u^{2\s+1}\,(u^2+1)^{\frac12(n-2)}\,du.\tag3.25$$
Note that the signature of $x^\perp$ is $(n-1,2)$, so that the volume
$\vol(\Gamma_x\back G_x)$ is also always finite. This proves Proposition .
\qed\qed\enddemo

\demo{\bf Proof of Proposition 3.4}

Finally, we consider the integral
$$
\int_{\Gamma_\ell\back D^+} 
\sum_{x\in \ell(\Q), \ x\ne0}
\ph(h^{-1}x)\,\b_{\s+1}(2\pi v R(x,z))\,
\,d\mu(z)\tag3.26$$
In this case, we choose a basis for $V$ such that the matrix for the 
inner product is
$$\pmatrix {}&{}&1\\{}&I_{n-1,1}&{}\\
1&{}&{}\endpmatrix\tag3.27$$
and such that $\ell$ is spanned by the first basis vector. Moreover, we assume that 
$$\{\ x\in \ell(\Q)\mid \ph(h^{-1}x)\ne0\ \}\subset 2\Z v_1.\tag3.28$$
The parabolic subgroup $P_\ell$ stabilizing the line $\ell$ then 
has Levi decomposition
$P_1 = U_1MA$ with $A\simeq GL(\ell)$, $M\simeq SO^+(n-1,1)$, 
and unipotent radical $U_1$. We take $z_0$ to be the oriented negative 
2--plane spanned by $\frac12(v_1-v_{n+2})$ and $v_{n+1}$ and let $K$ be its stabilizer. Then
$$G= SO^+(V)(\R) = UMAK\tag3.29$$ 
and we have the integral formula
$$\int_G \phi(g)\,dg = \int_U\int_M\int_A\int_K \phi(uma_rk)\,r^{-n-1}\,du\,dm\,dr\,dk,\tag3.30$$
where $a_r v_1 = r v_1$. For $z= uma_r\cdot z_0$, and $x = 2\a v_1$,  
$$R(x,z) = R(a_r^{-1}x,z_0) = 2 \a^2 r^{-2},\tag3.31$$
since 
$$a_r^{-1}x = 2 \a r^{-1} v_1 = 2\a r^{-1}\bigg( \frac12(v_1+v_{n+2}) + \frac12(v_1 - v_{n+2})\bigg)\tag3.32$$
has $\a\, r^{-1}\,(v_1 - v_{n+2})$ as its $z_0$ component. 
The integral (3.26) is then majorized by a constant times
$$\align
&\int_0^{\infty} 
\sum_{\a\in \Z, \ \a\ne0}
\,\b_{\s+1}(4\pi \a^2\, r^{-2} )\, r^{-n-1}\,dr\tag3.33\\
\nass
\nass
{}&=2(4\pi)^{-n}\,\zeta(n)\,\int_0^{\infty} 
\,\b_{\s+1}(r^{2} )\, r^{n-1}\,dr.
\endalign
$$
The integral here is finite provided $2\s+n\ge0$, so we obtain the required convergence 
provided $n\ge2$ and $\Re(\s) > -\frac{n}2$, i.e., in all isotropic cases except $n=1$ 
(which was an exceptional case). 
\qed\enddemo

\demo{\bf Proof of Proposition 3.6}
For $x$ with $Q(x)=m<0$, we write $x = \pr_z(x)+x'$ so that
$$R(x,z) = (x',x') - 2m \ge 2|m|.\tag3.34$$
We let $R'(x,z) = (x',x')$, and note that $R'(x,z)=0$ if and only 
if $x\in z$. 
Then we have the easy estimate:
$$\align
\int_1^\infty  e^{-2\pi v R(x,z)}\,v^{-\s-1}\,dv &\le e^{-2\pi R'(x,z)}\int_1^\infty  e^{-4\pi|m| v}\,v^{-\s-1}\,dv\\
\nass
{} &\le e^{-2\pi R'(x,z)}\int_1^\infty e^{-\e v}\, e^{(\e-4\pi|m|)v}\,v^{-\s-1}\,dv\tag3.35\\
\nass
{} &\le e^{-2\pi R'(x,z)}\,e^{\e-4\pi|m|}\,\int_1^\infty e^{-\e v}\, v^{-\s-1}\,dv\\
\nass
{} &\le C(\e,\s)\,e^{-2\pi R'(x,z)}\,e^{-4\pi|m|},
\endalign
$$
for any $\e$ with $0<\e<4\pi|m|$, where the constant $C(\e,\s)$ is uniform in any $\s$--halfplane and 
independent of $m$. Note that there is a positive lower bound for the quantity $|m|$ where $m<0$ has $c_\ph(-m)\ne0$. 
This leads to the expression
$$\sum_{m<0} |c_\ph(-m)|\,e^{-4\pi|m|}\,
\sum_{\scr x\atop {\scr Q(x)=m \atop \scr \mod \Gamma}} \ph(h^{-1}x)\,
\int_{\Gamma_x\back D^+}e^{-2\pi R'(x,z)}\,d\mu(z).\tag3.36$$

Recall that the modular form $f_\ph$ with Fourier coefficients $c_\ph(-m)$ has weight $1-\frac{n}2$, 
with some real multiplier for a congruence subgroup of $\text{\rm SL}_2(\Z)$,  
and is holomorphic in the upper halfplane with possible poles at the cusps. 
Then it is known that 
$$c_\ph(-m) = O\bigg(\  |m|^{-\frac{n+1}4}\, e^{C \sqrt{|m|}}\ \bigg),\tag3.38$$
i.e., these coefficients grow at most subexponentially. 
The (explicit) constant $C$ depends only on the order of the pole of $f_\ph$ and on the 
multiplier.  
If $n>2$, so that 
$f_\ph$ has negative weight, this fact follows from the classical work of Rademacher \cite{\rademacher}, 
Rademacher--Zuckermann \cite{\rademacherZ},
Zuckermann \cite{\zuckerman}, and Petersson\cite{\petersson}, 
cf. also Hejhal \cite{\hejhal}. The cases $n=1$ and $2$ are covered by Hejhal \cite{\hejhal} and Niebur \cite{\niebur}. 

Finally, it remains to estimate the quantity
$$\sum_{\scr x\atop {\scr Q(x)=m \atop \scr \mod \Gamma}} \ph(h^{-1}x)\,
\int_{\Gamma_x\back D^+}e^{-2\pi R'(x,z)}\,d\mu(z).\tag3.39$$
To estimate the integral here, 
we choose basis for $V$ so that the inner product has matrix
$-I_{2,n}$ and such that $x=2\a v_1$. Let $z_0$ be the span of $v_1$ and $v_2$, 
and let $A=\{a_t\}$ be the 1--parameter subgroup which is the identity 
component of the special orthogonal group of the plane spanned by $v_1$ and $v_3$.
In this case $a_t v_1 =\cosh(t)v_1+\sinh(t)v_3$. Again we have the 
decomposition (3.29) and an integral formula analogous to (3.30), but with the $\cosh$ and 
$\sinh$ switched in the modulus factor.  For $z=g_xa_t\cdot z_0$, we also have
$$R(x,z) = 2|m|\cosh^2(t).\tag3.40$$
and
$$R'(x,z) = 2|m|\cosh^2(t) - 2|m| = 2|m|\sinh^2(t).\tag3.41$$
Then we have
$$\align
&\int_{\Gamma_x\back D^+}  e^{-2\pi R'(x,z)}\,d\mu(z)\\
\nass
{}&=C'\, \vol(\Gamma_x\back G_x)\,\vol(K)
\int_0^\infty e^{-4\pi |m|\,\sinh^2(t)}\,\sinh(t)^{n-1}\cosh(t)\,dt\tag3.42\\
\nass
\nass
{}&=C'\, \vol(\Gamma_x\back G_x)\,\vol(K)\,
\frac12\, (4\pi |m|)^{-\frac{n}2}\,\Gamma(\frac{n}2).
\endalign
$$ 
Using this in (3.39), we are left with
$$\sum_{m<0} |m|^{-\frac{3n+1}{4}}\,e^{C \sqrt{|m|}-4\pi|m|}
\sum_{\scr x\atop {\scr Q(x)=m \atop \scr \mod \Gamma}} \ph(h^{-1}x)\,
\vol(\Gamma_x\back G_x).\tag3.43
$$
Here $m$ runs over the negative elements of $N^{-1}\Z$ for a suitable $N$
depending on $\ph$ and $h$. The resulting expression is finite
since, \cite{\siegeltata},
$$\sum_{\scr x\atop {\scr Q(x)=m \atop \scr \mod \Gamma}} \ph(h^{-1}x)\,
\vol(\Gamma_x\back G_x)  = O(|m|^{\frac{n}2+\e}).\tag3.44
$$
\qed\enddemo

This completes the proof of Theorem~3.1. 
\qed\qed\qed

There are several more things which need to be proved. 

\demo{Proof of Proposition~2.5}  By (2.3), the left hand side of the identity 
$$\multline
\CT{\s=0}\bigg\{\lim\limits_{T\rightarrow\infty} \int_{\Cal F_T}\sum_\ph f_\ph(\tau)\,I(\tau,\ph^\ev)\, v^{-\s-2}\,du\,dv\bigg\}\\
\nass
{}= \lim\limits_{T\rightarrow \infty}\bigg[\  \int_{\Cal F_T}\sum_\ph f_\ph(\tau)\,I(\tau,\ph^\ev)\, v^{-2}\,du\,dv
-c_0(0)\,\log(T)\ \bigg].
\endmultline
$$ 
to be proved can be written as
$$\align
&\vol(X)^{-1}\,\CT{\s=0}\bigg\{\lim\limits_{T\rightarrow\infty} \int_{X_K}\int_{\Cal F_T} ((\,f(\tau),\vartheta(\tau,z)\,))\, v^{-\s-2}\,du\,dv\,
d\mu(z)\bigg\}\\
\nass
{}&= \vol(X)^{-1}\,\int_{X_K}\int_{\Cal F_1} ((\,f(\tau),\vartheta(\tau,z)\,))\, v^{-2}\,du\,dv\,
d\mu(z)\\
\nass
{}&\qquad\qquad+\vol(X)^{-1}\,\CT{\s=0}\bigg\{\lim\limits_{T\rightarrow\infty} \int_{X_K} \int_1^T C(v,z)\,v^{-\s-1}\,dv\,d\mu(z)\,\bigg\}
\endalign$$
The analysis made in the proof of Theorem~3.1 above shows that the integral
$$\int_{X_K} \int_1^\infty \bigg[\,C(v,z)- C_{00}(v,z)\,\bigg]\,v^{-\s-2}\,dv\,d\mu(z)$$
defines a holomorphic function of $\s$ in the half plane $\Re(\s)>-1$. Note that in the case $n=1$
there are no $C_0$ terms, since $V$ is then assumed to be anisotropic. The remaining term 
is
$$\vol(X)^{-1}\int_{X_K}\int_1^T C_{00}(v,z)\,v^{-\s-1}\,dv\,d\mu(z) = c_0(0)\,\frac{1}{\s}\big(1-T^{-\s}\big) = c_0(0)\log(T) + O(\s).$$
This term makes no contribution when we take the limit as $T$ goes to infinity followed by the 
constant term at $\s=0$. Thus, once the term $c_0(0)\,\log(T)$ has been removed, we can pass to the limit on $T$ with $\s=0$, and this
proves Proposition~2.5.
\qed\enddemo

\demo{Proof of Proposition~2.11} In the Fourier expansion (2.21) for 
$E(\tau,s;\ph,\ell+2)$ for a factorizable function $\ph=\tt_p\ph_p\in S(V(\A_f))$, the 
$m$th coefficient, for $m\ne0$, has a product formula
$$E_m(\tau,s;\ph,\ell+2) = A_\ph(s,m,v)\,q^m = W_{m,\infty}(\tau,s;\ell+2)\cdot\prod_{p}W_{m,p}(s,\ph_p).$$
The following facts are well known, cf. \cite{\ky} for more details. 
For $s=s_0=\ell+1 = \frac{n}2$, 
$$\align
W_{m,\infty}(\tau,\frac{n}2;\frac{n}2+1) &= \frac{(-2i)^{\frac{n}2+1}}{\Gamma(\frac{n}2+1)}\, m^{\frac{n}2}\,q^m\\
\noalign{if $m>0$,}
W_{m,\infty}(\tau,\frac{n}2;\frac{n}2+1) &=0\\
\noalign{ if $m<0$, and}\\
W'_{m,\infty}(\tau,\frac{n}2;\frac{n}2+1) &= \pi (-i)^{-\frac{n}2-1}\,2^{-\frac{n}2}\, 
q^m\,v^{-\frac{n}2}\,\int_1^\infty e^{-4\pi|m|vr}\,r^{-\frac{n}2-1}\,dr.
\endalign$$ 
On the other hand, for any $m\ne0$, the product over the finite primes is 
$$C(m):=\bigg(\,\prod_{p}W_{m,p}(s,\ph_p)\,\bigg)_{s=s_0} = O(1).$$
Therefore, for $m<0$, we have
$$b_\ph(m,v) = \pi (-i)^{-\frac{n}2-1}\,2^{-\frac{n}2}\,C(m)\,v^{-\frac{n}2}\,\int_1^\infty e^{-4\pi|m|vr}\,r^{-\frac{n}2-1}\,dr,$$
where $C(m) = O(1)$. 
Thus
$$|b_\ph(m,v)| = O\big(v^{-\frac{n}2-1}\,|m|^{-1} \,e^{-4\pi|m|v}\,\big).$$
Using (3.38), this proves part (i) and (ii) of Proposition~2.11.

Finally, the constant term has the form
$$E_0(\tau,s;\ph,\ell+2) = v^{\frac12(s+1-\ell)}\,\ph(0) + W_{0,\infty}(\tau,s;\ell+2)\,\prod_{p} W_{0,p}(s,\ph_p),$$
where
$$W_{0,\infty}(\tau,s;\ell+2) = 2\pi\,(-i)^{\frac{n}2+1}\,v^{-\frac12(s+\frac{n}2)}\,2^{-s}\frac{\Gamma(s)\frac12(s-\frac{n}2)}{\Gamma(\frac12(s+\frac{n}2+2))
\Gamma(\frac12(s-\frac{n}2+2))}.$$
Then, the derivative at $s=s_0=\frac{n}2$ is
$$E'_0(\tau,\frac{n}2;\ph,\frac{n}2+1) = \frac12\,\log(v)\,\ph(0) +
\pi\,(-i)^{\frac{n}2+1}\,v^{-\frac{n}2}\,2^{-\frac{n}2}\frac{\Gamma(\frac{n}2)}{\Gamma(\frac{n}2+1)}\,C(0).$$
This yields (iii) of Proposition~2.11.
\qed\enddemo

\subheading{\Sec4. Formulas for degrees}

In this section, we explain how the Siegel--Weil formula can 
be applied to yield formulas for the degrees of certain divisors on the 
quasiprojective varieties attached to orthogonal groups of signature $(n,2)$ over $\Q$. 
More precisely, these degrees occur as the Fourier coefficients of certain 
(special values of) Eisenstein series. The basic idea will be to apply the 
Siegel--Weil formula for two different quadratic spaces to describe a special value of 
the same Eisenstein series! Comparison of the Fourier coefficients of the two theta 
integrals and the Eisenstein series yields nontrivial identities, several of which 
occur in the classical literature, \cite{\zagier}, \cite{\cohen}, \cite{\vdgeer}. 

\subheading{The Siegel-Weil formula}

For convenience of the reader, we briefly review the Siegel--Weil formula for 
the dual pair $(SL_2,O(V))$ needed in this section and in section 2. 

Let $V$ be a nondegenerate quadratic space over $\Q$, and 
let $G=SL_2$. As before, let $G'_\A$ be the metaplectic cover of $G(\A) = SL_2(\A)$. 
We identify $G'_\A = SL_2(\A)\times \{\pm1\}$, where multiplication on the right 
is given by $[g_1,\e_1][g_2,\e_2] = [g_1g_2,\e_1\e_2 c(g_1,g_2)]$, for the cocycle 
as in \cite{\waldspurger}, \cite{\gelbart}. In particular, we have subgroups
$$N_\A = \{ n=[n(b),1]\mid b\in \A\},\qquad n(b) = \pmatrix 1&b\\{}&1\endpmatrix\tag4.1$$
and 
$$M_\A = \{ \und{m}(a)= [m(a),\e]\mid a\in\A^\times, \e=\pm1\},\qquad m(a) = \pmatrix a&{}\\{}&a^{-1}\endpmatrix.\tag4.2$$
An idele character $\chi$ of $\Q^\times\back \A^\times$ determines a character $\chi^\psi$ of 
$M_\A$ by
$$\chi^\psi([m(a),\e]) = \e\,\chi(a)\,\gamma(a,\psi)^{-1}\tag4.3$$
where $\gamma(\cdot,\psi)$ is the global Weil index.

The group $G'_\A$ acts on the Schwartz space $S(V(\A))$ via the Weil representation $\o=\o_\psi$
determined by our fixed additive character $\psi$ of $\A/\Q$, and 
this action commutes with the linear action of $O(V)(\A)$. For $g'\in G'_\A$,  
$h\in O(V)(\A)$, and $\ph\in S(V(\A))$, the theta series
$$\theta(g',h;\ph) = \sum_{x\in V(\Q)} \o(g')\ph(h^{-1}x),\tag4.4$$
defines a smooth function on $G'_\A\times O(V)(\A)$, left invariant under $G'_\Q\times O(V)(\Q)$, 
and slowly increasing on the quotient $\big(G'_\Q\times O(V)(\Q)\big) \back \big(G'_\A\times O(V)(\A)\big)$. 

By Weil's criterion \cite{\annals} in the present case, the theta integral 
$$I(g';\ph) = \int_{O(V)(\Q)\back O(V)(\A)} \theta(g',h;\ph)\, dh,\tag4.5$$
where $\vol(O(V)(\Q)\back O(V)(\A),dh)=1$, is absolutely convergent whenever either $V$ is 
anisotropic or $\dim(V)-r >2$, where $r$ is the Witt index of $V$. 
The resulting automorphic form $I(\ph)$ on $G'_\Q\back G'_\A$ is identified, by the 
Siegel--Weil formula, with a special value of an Eisenstein series, defined as follows. 

Let $\chi=\chi_V$ be the quadratic character of $\A^\times/\Q^\times$ defined by
$$\chi(x) = (x,(-1)^{m(m-1)/2}\det(V)),\tag4.6$$
where $m=\dim(V)$ and $\det(V)\in \Q^\times/\Q^{\times,2}$ is the determinant of the 
matrix for the quadratic form $Q$ on $V$. For $s\in \C$, let 
$I(s,\chi)$ be the principal series representation of $G'_\A$ 
consisting of smooth functions $\P(s)$ on $G'_\A$ 
such that 
$$\P(n\, \und{m}(a) g',s) = \cases \chi^\psi(\und{m}(a))\,|a|^{s+1}\,\P(g',s)&\text{ if $n$ is odd,}\\
\nass
\\ \chi(m(a))\,|a|^{s+1}\,\P(g',s)&\text{ if $n$ is even.} 
\endcases\tag4.7
$$ 
There is then a $G'_\A$ intertwining map
$$\l=\l_V:S(V(\A)) \lra I(s_0,\chi_V),\qquad \l(\ph)(g') = \o(g')\ph(0),\tag4.8$$
where $s_0 = \frac{m}2-1$. As in section 1, let $K'_\infty K'$ be the full inverse image 
of $SO(2)\times \SL_2(\hat\Z)$ in $G'_\A$. A section $\P(s)\in I(s,\chi)$ will be called 
standard if its restriction to $K'_\infty K'$ is independent of $s$. By the Iwasawa decomposition
$G'_\A= N'_\A M'_\A K'_\infty K'$, the function $\l(\ph)\in I(s_0,\chi)$ 
has a unique extension to a standard section $\P(s)\in I(s,\chi)$, where $\P(s_0) =\l(\ph)$.
The Eisenstein series, defined by
$$E(g',s;\P) = E(g',s;\ph) = \sum_{\gamma\in P'_\Q\back G'_\Q} \P(\gamma g',s)\tag4.9$$ 
for $\Re(s)>1$ has a meromorphic analytic continuation to the whole $s$--plane. 
\proclaim{Theorem 4.1} {\rm (Siegel--Weil formula)} (i) Assume that 
$V$ is anisotropic or that $\dim(V)-r >2$, where $r$ is the Witt index of $V$, 
so that the theta integral (4.5) is absolutely convergent. Then $E(g',s;\ph)$ is 
holomorphic at the point $s=s_0 = \frac{m}2-1$, where $m=\dim(V)$, and 
$$E(g',s_0;\ph) = \kappa\cdot I(g',\;\ph),$$
where $\kappa=2$ when $m\le2$ and $\kappa=1$ otherwise. 
\hfill\break
(ii) Suppose, in addition, that $m>1$. Then
$$E(g',s_0;\ph)=\kappa\cdot I(g';\ph) = \frac{\kappa}2\int_{SO(V)(\Q)\back SO(V)(\A)} \theta(g',h;\ph)\,dh,$$
where $dh$ is Tamagawa measure on $SO(V)(\A)$. 
\endproclaim

When $m>4$ part (i) is the classic result of Siegel and Weil, in Weil's formulation. The variants 
for $m\le4$ are also mostly classical, e.g., due to Hecke, Siegel, etc.. We do not attempt to 
give systematic references. Of course, the analogous result holds for any number field $F$.

The point of (ii) is that, we can almost always replace the integral over $O(V)(\Q)\back O(V)(\A)$ 
with the integral over $SO(V)(\Q)\back SO(V)(\A)$. The later is much more convenient, since $SO(V)$ 
is connected. In the range $m>4$, this fact is again a very special case of the results 
of Weil, \cite{\weilII}, pp.76--77, Th\'eor\`em 5. 

\demo{Proof} We only prove part (ii). Let $I':S(V(\A))\rightarrow \C$
be the linear functional given by  
$$I'(\ph) = \int_{SO(V)(\Q)\back SO(V)(\A)} \theta(h;\ph)\,dh,\tag4.10$$
where $\theta(h;\ph) = \theta(e,h;\ph)$, so that $I'$ defines an element of 
$$\text{\rm Hom}_{SO(V)(\A)}\big(\,S(V(\A)),\C\,\big),\tag4.11$$
where $SO(V)(\A)$ acts trivially on $\C$. 
The group 
$$C(\A_f) = O(V)(\A)/SO(V)(\A) \simeq \mu_2(\A)\tag4.12$$ 
acts on this space of such functionals. In fact one has:
\proclaim{Proposition 4.2} If $\dim(V)>1$, then the action of $C(\A_f)$ on the space of 
$SO(V)(\A)$--invariant linear functionals on $S(V(\A))$
(4.11) is trivial.
\endproclaim
\demo{Proof} For any prime $p\le\infty$, consider the analogous local space
$$Hom_{SO(V_p)}(S(V_p),\C),\tag4.13$$
with its action of $C_p$. If the sign character $\e_p$ of 
$C_p=O(V_p)/SO(V_p)$ occurs, then the sign representation $\sgn_p$ of $O(V_p)$ 
occurs in the local theta correspondence for the dual pair 
$(\widetilde{\SL}_2(\Q_p),O(V_p))$. But it is known, \cite{\rallisHDC}, that the sign representation
does not occur for such a dual pair if $\dim(V)=m>1$. Thus $C_p$ acts trivially on (4.13), and a standard 
argument then shows that $C(\A)$ acts trivially on (4.11), as claimed. 
\qed\enddemo
On the other hand, it is clear that
$$\align
I(g',\ph) &= \int_{C(\Q)\back C(\A)} I'\big(\o(g')\o(h)\ph\big)\,dc\\
\nass
{}&= \frac12\, \int_{C(\A)}I'\big(\o(g')\o(h)\ph\big)\,dc\tag4.14\\
\nass
{}&= \frac12\,I'\big(\o(g')\ph\big),
\endalign
$$
where $h\in O(V)(\A)$ projects to $c\in C(\A)$ and where $\vol(C(\A),dc)=1$. The factor $1/2$ occurs
as the volume of $C(\Q)\back C(\A)$. 
\qed\qed\enddemo

\subheading{A matching principle} 

For a nondegenerate quadratic space $V$ over $\Q$ of dimension $m$, let
$$\Pi(V) = \text{ image}(\lambda_V)\subset I(s_0,\chi)\tag4.15$$
be the resulting $G'_\A$-submodule of the principal series, where $\chi=\chi_V$ and $s_0=\frac{m}2-1$. There are analogous local maps
$$\l_p:S(V_p) \lra I_p(s_0,\chi_p),\tag4.16$$
with images
$$\Pi_p(V_p) = \text{ image}(\lambda_p)\subset I_p(s_0,\chi_p),\tag4.17$$
the local components of $\Pi(V)$ for the corresponding local induced representations. Note that $\Pi(V)$ and the $\Pi_p(V_p)$'s 
are not always irreducible. The key idea is that the Eisenstein series (4.9) associated to 
$\ph=\tt_p\ph_p\in S(V(\A))$ depends only on the the collection $\{\l_p(\ph_p)\}$ 
of local components. 

{\bf Definition 4.3:} Let $V_p$ and $V'_p$ be quadratic spaces over $\Q_p$ of dimension $m$ 
and fixed character $\chi_{V_p}=\chi_{V'_p}=\chi_p$. Functions $\ph_p\in S(V_p)$ and 
$\ph'_p\in S(V'_p)$ are said to {\bf match} if
$$\l_p(\ph_p) = \l'_p(\ph'_p).$$

{\bf Remark  4.4:} This matching is analogous to that which occurs in the trace formula 
and relative trace formula, and our identity of theta integrals can be viewed as an 
analogue of a comparison of trace formulas. 

{\bf Examples 4.5:}
If $m>4$, or if $m=4$ and $\chi_p\ne 1$, then the nonarchimedean local principal series 
$I_p(s_0,\chi_p)$ are irreducible and hence, for any pair $V_p$ and $V'_p$, 
every $\ph_p\in S(V_p)$ has a matching $\ph'_p\in S(V'_p)$. 

If $m=4$, and $\chi_p = 1$, then $s_0=1$ and $I_p(s_0,\chi_p)$ has the special representation 
as irreducible submodule and the trivial representation as quotient. The split $4$ dimensional 
quadratic space $V_p$ has $\Pi_p(V_p) = I_p(s_0,\chi_p)$, while the anisotropic space $V'_p$ 
given by the reduced norm on the division quaternion algebra over $\Q_p$ has $\Pi_p(V'_p)$ 
the irreducible special. Therefore the space of $\ph_p$'s in $S(V_p)$ which have 
matching $\ph'_p$'s has codimension $1$. 

If $m=3$, then $s_0=\frac12$ and $I_p(s_0,\chi_p)$ always has length $2$ with a special representation of $G'_p$ as the 
irreducible subrepresentation and an irreducible Weil representation -- playing the role of 
the trivial representation for the metaplectic group $G'_p$ -- as irreducible quotient, \cite{\rallisschiff}. 
The ternary quadratic space $V_p$ of trace $0$ elements in $M_2(\Q_p)$ with a scalar multiple 
(determined by $\chi_p$) of 
the determinant form has $\Pi_p(V_p)=I_p(s_0,\chi_p)$, and the analogous space $V'_p$ of trace $0$ elements 
in the division quaternion algebra over $\Q_p$ has $\Pi_p(V'_p)\subset I_p(s_0,\chi_p)$ 
the unique irreducible subrepresentation. Now the subspace of $\ph_p$'s in $S(V_p)$ which have 
matching $\ph'_p$'s in $S(V'_p)$ has infinite codimension. 

{\bf Remark 4.6:} If $m=2$ and $\chi_p\ne1$, then the spaces $\Pi_p(V_p)$ and $\Pi_p(V'_p)$ are irreducible and 
distinct, while, if $m=1$, there is a unique space with a given $\chi_p$, so the
matching phenomenon of interest here will not occur globally. 

Note that the cases $m=3$ and $4$ are precisely those for which $s_0$ is in or at the edge of the 
critical strip $|\Re(s)|\le 1$.  

Over $\R$, the situation is the following. For 
$r\in \frac12\Z$, satisfing a 
suitable parity condition,
let $\P^r(s)$ be the (unique) function in $I_\infty(s,\chi_\infty)$ such that 
$$\P^r(k',s) = \chi_r(k'),\tag4.18$$
for the character $\chi_r$ of $K'_\infty$. The space of $K'_\infty$--finite vectors in $I_\infty(s,\chi_\infty)$ 
is then spanned by the $\P^r(s)$'s for $r\in r_0+ 2\Z$.  
\proclaim{Lemma 4.7}  Suppose that $V_\infty$ and $V'_\infty$ are quadratic spaces over $\R$ 
of dimension $m$ and with the same quadratic character, i.e., with signatures $(p,q)$ and $(p',q')$ 
with $q\equiv q'\mod(2)$. Suppose that $\ph_\infty\in S(V_\infty)$ and $\ph_\infty'\in S(V_\infty')$ 
are eigenfunctions for $K'_\infty$ with eigencharacter $\chi_r$ and with
$\ph_\infty(0) = \ph_\infty'(0)$. Then $\ph_\infty$ and $\ph'_\infty$ 
match and $\lambda_\infty(\ph_\infty) = \lambda_\infty'(\ph_\infty') = \P^r(s_0)$. 
\endproclaim

\proclaim{Matching Principle 4.8} Suppose that $V$ and $V'$ are quadratic spaces over $\Q$ of the same dimension and 
with the same quadratic character $\chi_V=\chi_{V'}=\chi$. Suppose that $\ph\in S(V(\A))$ and 
$\ph'\in S(V'(\A))$ match, i.e., $\lambda(\ph)=\lambda'(\ph')=\P(s_0)$.  Assume that the 
convergence condition of the Siegel--Weil 
formula is fulfilled by the spaces $V$ and $V'$. Then
$$I(g',\ph) = E(g',s_0,\P)=I(g',\ph').$$
\endproclaim

{\bf Remark 4.9.}  The definition of matching and the resulting 
equality of theta integrals can be extended to dual pairs $(Sp(r),O(V))$, $(Sp(r),O(V'))$ for any $r\ge1$ 
over a number field, dual pairs for unitary groups etc.,  
etc.  

Of course, the matching principle is a trivial observation, but, while the Eisenstein 
series is built from purely local date, the theta integrals involved depend on global arithmetic. 
In particular, their equality can yield some highly nontrivial identities. We now describe one of these. 

\subheading{ A geometric example}

Let $V$ be a quadratic 
space over $\Q$ of signature $(n,2)$, and let $V'$ be a quadratic space over $\Q$ 
with $\chi_{V'}=\chi_V=\chi$ but with signature $(n+2,0)$.  
Suppose that $\ph\in S(V(\A_f))$ and $\ph'\in S(V'(\A_f))$ are matching functions. 
By the discussion above, when $n>2$, any $\ph$ has such a matching $\ph'$. We 
next construct matching functions over $\R$.

As explained in section 1 above, the Gaussian for $V(\R)$ is  
the function $\ph_\infty\in S(V(\R))\tt A^{(0,0)}(D)$ given by 
$$\ph_\infty(x,z) = e^{-\pi(x,x)_z} = e^{-2\pi R(x,z)}\,e^{-2\pi Q(x)}.\tag4.19$$
It has weight $\ell=\frac{n}2-1$ and $\ph_\infty(0,z)=1$, so that 
$$\lambda_\infty(\ph_\infty) = \P^\ell(s_0).\tag4.20$$

Let $V'(\R)$ be a quadratic space of signature $(n+2,0)$. 
The Gaussian $\ph_\infty'\in S(V'(\R))$ is given by
$$\ph'_\infty(x) = e^{-2\pi Q'(x)}.\tag4.21$$
It has weight $\frac{n+2}2=\ell+2$ and $\ph'_\infty(0)=1$, so that 
$$\lambda'_\infty(\ph_\infty') = \P^{\ell+2}(s_0).\tag4.22$$
In particular, the Gaussians of $V(\R)$ and $V'(\R)$ {\it do not match}, and we will 
need to find another function for $V(\R)$. 

One of the main results of \cite{\kmillsonI} was the construction of a 
Schwartz {\it form} for $V$, 
$$\ph_{KM}\in S(V(\R))\tt A^{(1,1)}(D),\tag4.23$$
where $A^{(1,1)}(D)$ is the space of smooth $(1,1)$--forms on $D$, with the following properties:
\roster
\item"{(i)}"  For all $h\in O(V(\R))$, 
$$h^*\ph_{KM}(h^{-1}x) = \ph_{KM}(x),\tag4.24$$ 
where $h^*$ indicated the action of $h$ on the space $A^{(1,1)}(D)$ by pullback. 
\item"{(ii)}"  $\ph_{KM}$ has weight $\ell+2$ for $K'_\infty$, i.e., 
$$\o(k')\ph_{KM} = \chi_{\ell+2}(k')\,\ph_{KM},\tag4.25$$
for the Weil representation action of $K'$ on $S(V(\R))$. 
\item"{(iii)}"
$\ph_{KM}$ is closed:
$$d\ph_{KM}=0\tag4.26$$
for exterior differentiation $d$ on $D$. 
\endroster

Note that it follows from properties (i) and (iii) above that $\ph_{KM}(x)\in A^{(1,1)}(D)$ 
is a closed $O(V(\R))_x$ invariant form. For example, 
$$\O:=\ph_{KM}(0)\tag4.27$$
is an $O(V(\R))$ invariant $(1,1)$--form on $D$, which 
we will identify in a moment.

In the present situation, $\ph_{KM}$ is obtained as follows. Recall from Lemma~3.8 that for $t\in \R_{>0}$, 
the exponential integral $\b_1(t)$ has a logarithmic 
singularity, $-\log(t)$,  as $t\rightarrow 0$ and decays like $e^{-t}$ as $t\rightarrow \infty$. 
For $x\in V(\R)$, $x\ne 0$, and $z\in D$, let 
$$\xi(x,z) =\b_1(2\pi R(x,z))\,  e^{-2\pi Q(x)}.\tag4.28$$
This function is smooth away from the incidence locus
$$\{[x,z] \in V(\R)\times D\mid \pr_z(x) =0\,\}.\tag4.29$$
For example, if $x\in V(\R)$ is fixed, then $\xi(x)$ is a 
smooth function on $D-D_x$, 
where
$$D_x=\{\ z\in D\mid z\perp x\,\},\tag4.30$$
as in (1.44). 
Moreover, $\xi(x,z)$ decays exponentially as $z$ goes to infinity away from $D_x$. 
Note that $D_x$ is nonempty if and only if $Q(x)>0$. 
The crucial fact then is that,  for $x\ne 0$, 
$$\ph_{KM}(x) = dd^c \xi(x),\tag4.31$$
where, $d^c= \frac{1}{4\pi i}(\partial-\bar\partial)$. 
In fact, as in \cite{\annals}, we have the stronger assertion:
\proclaim{Proposition 4.10} As currents on $D$, 
$$d d^c\xi(x) + e^{-2\pi Q(x)}\,\delta_{D_x} = [\ph_{KM}(x)].$$
\endproclaim
\demo{Proof} Omitted
\qed\enddemo

We can recover the explicit formula for $\O$ from this result. 
\proclaim{Proposition 4.11} On the tube domain $\Bbb D$, let $\rho=\rho(z) = -\frac12(w(z),w(\bar z))$, 
be the norm of the section $z\mapsto w(z)$ of $\Cal L_D$, as in (1.10). Then
$$\align
\O &= dd^c\log(\rho)\\
\nass
{}&=-\frac{1}{2\pi i}\bigg[\ -(y,y)^{-2}(y,dz)\wedge (y,d\bar z) + (y,y)^{-1}\frac12(dz,d\bar z)\ \bigg].
\endalign
$$
\endproclaim
\demo{Proof} We compute
$$\align
dd^c\xi(x) &= -\frac{1}{2\pi i} \d\db\bigg\{ \b_1(2\pi R)\bigg\}\cdot e^{-2\pi Q(x)}\\
\nass
{}&= \frac{1}{2\pi i} \d\bigg\{ e^{-2\pi R} \db \log(R)\bigg\}\cdot e^{-2\pi Q(x)}\tag4.32\\
\nass
{}&= \frac{1}{2\pi i}\bigg[ -2\pi \d R\wedge\db \log(R) + \d\db \log(R)\bigg]\cdot e^{-2\pi R-2\pi Q(x)}\\
\nass
{}&=\ph_{KM}(x)
\endalign
$$
For a moment, we write $\a=(x,w(z))$ and $\rho=|y|^2=-(y,y)$, as in (1.10), so that, by (1.16), 
$R= \rho^{-1}|\a|^2$. Then
$$\d\db\log(R) = - \d\db\log(\rho),\tag4.33$$
and
$$\d R\wedge \db\log(R) = \rho^{-1}d\a\wedge d\bar\a
 -\rho^{-2}\,\bar\a\,d\a\wedge \db\rho -\rho^{-2}\,\a\,\d\rho\wedge d\bar\a +\rho^{-3}|\a|^2\d\rho\wedge\db\rho.\tag4.34$$
Notice that this last expression defines a smooth form on $D$. 

Setting $x=0$, we obtain:
$$\O = \ph_{KM}(0) = dd^c\log(\rho).\tag4.35$$
But now, writing
$$\rho = -(y,y) = \frac14(z-\bar z,z-\bar z),\tag4.36$$
we have
$$\align
\O&= -\frac{1}{2\pi i} \d\db \log(\rho)\\
\nass
{}&= -\frac{1}{2\pi i}\bigg[\ -\rho^{-2}\,\d\rho\wedge\db\rho + \rho^{-1}\d\db\rho\ \bigg]\tag4.37\\
\nass
{}&= -\frac{1}{2\pi i}\bigg[\ -(y,y)^{-2}(y,dz)\wedge (y,d\bar z) + (y,y)^{-1}\frac12(dz,d\bar z)\ \bigg]
\endalign
$$
as claimed.\qed\enddemo

\proclaim{Corollary 4.12} The form 
$$\O =\ph_{KM}(0)= dd^c\log||s||^2$$
on $X_K$ is the  first Chern form for the 
holomorphic line bundle $\Cal L^\vee$ dual to $\Cal L$. 
In particular, $-\O$ is an invariant K\"ahler form on $\Bbb D$ and hence determines 
a K\"ahler form on $X_K$. 
\endproclaim

{\bf Examples 4.13:} In the case $n=1$ we have $\Bbb D \simeq \C\setminus \R = \H^+\cup \H^-$ and 
$\O= -\frac1{2\pi}\, y^{-2}\,dx\wedge dy$. In the case $n=2$, we have 
$\Bbb D \simeq \H\times\H$ and 
$$\O = -\frac1{4\pi}\bigg(\, y_1^{-2}\,dx_1\wedge dy_1
+y_2^{-2}\,dx_2\wedge dy_2\,\bigg),\tag4.38$$
(compare \cite{\hirzebruchzagier}, p.104, \cite{\vdgeerbook}, p.102.)

We now return to the theta integral and its geometric meaning. 
Write
$$\ph_{KM}(x)\wedge \O^{n-1} = \tilde{\ph}_{KM}(x)\,\O^{n},\tag4.39$$
for a function $\tilde{\ph}_{KM}\in S(V(\R))\tt A^{(0,0)}(D)$. Note that, 
since
$\O$ is $O(V(\R))$--invariant,
$$\tilde\ph_{KM}(hx,hz)=\tilde\ph_{KM}(x,z)\tag4.40$$
for all $h\in O(V(\R))$. Moreover,  $\tilde{\ph}_{KM}$ 
also has weight $\ell+2$ for the Weil representation action of $K'_\infty$. 

\proclaim{Lemma 4.14}  For all $z\in D$, 
$$\lambda(\ph_{KM}(\cdot, z)) = \P^{\ell+2}(s_0)\,\O,$$
and
$$\lambda(\tilde\ph_{KM}(\cdot,z)) = \P^{\ell+2}(s_0).$$
\endproclaim

\proclaim{Corollary 4.15} For all $z\in D$, the functions $\tilde{\ph}_{KM}(\cdot, z)\in S(V(\R))$ and 
$\ph_0'\in S(V'(\R))$ match. 
\endproclaim

We now return to the global situation, so that, for the matching functions 
$\ph\in S(V(\A_f))$ and $\ph'\in S(V'(\A_f))$ above,
$$\lambda(\tilde{\ph}_{KM}\tt\ph) = \lambda'(\ph'_0\tt\ph) =\P(s_0),\tag4.41$$
for a standard section $\P(s)\in I(s,\chi)$. 
Hence, we have an equality of Eisenstein series:
$$E(g',s,\lambda(\tilde{\ph}_{KM}\tt\ph)) = E(g',s,\lambda'(\ph'_0\tt\ph'))= E(g',s,\P).\tag4.42$$
Applying the Siegel--Weil formula, we have
\proclaim{Corollary 4.16} 
$$I(g',\tilde{\ph}_{KM}\tt\ph;V) =  I(g',\ph'_0\tt\ph';V') = E(g',s_0,\P).$$
\endproclaim
Here, in forming the theta integral of $V$, we use the theta function 
$$\theta(g',h;\tilde\ph_{KM}\tt\ph) =\sum_{x\in V(\Q)} \o(g')\tilde\ph_{KM}(h_\infty^{-1}x,z_0)\,\ph(h^{-1}x),\tag4.43$$
on $G'_\A\times O(V(\A))$, where $z_0\in D$ is a fixed point. In particular, as a function on 
$O(V(\A))$ this function is right invariant under the stabilizer of $z_0$ in $O(V(\R))$.

Next we would like to explain the geometric content of the first of these expressions. 
The key point is to determine the relation between the integral of the theta function (4.5) over 
$O(V)(\Q)\back O(V)(\A)$ and the integral of the differential form $\theta(g',\ph_{KM}\tt\ph)\wedge \O^{n-1}$ 
over $X_K$.

\proclaim{Proposition 4.17} Assume that the compact open subgroup $K\subset H(\A_f)$ satisfies:
$$Z_K:=K\cap Z(\A)  \simeq \hat\Z^\times,\tag Z$$
under the isomorphism $Z(\A)\simeq \A^\times$. Then 
$$
I(g',\tilde{\ph}_{KM}\tt\ph;V) =(-1)^n\frac14\,\vol(K) \int_{X_K}\theta(g',\ph_{KM}\tt\ph)\wedge \O^{n-1}.
$$
\endproclaim
\demo{Proof}
By (ii) of Theorem~4.1, 
$$I(g';\tilde\ph_{KM}\tt\ph;V) = \frac12\,\int_{SO(V)(\Q)\back SO(V)(\A)} \theta(g',h;\tilde\ph_{KM}\tt\ph)\,dh\tag4.44$$
where $dh$ is Tamagawa measure on $SO(V)(\A)$. 
A factorization $dh= dh_\infty\times dh_f$ will be determined by the choice of $dh'_\infty$ made below. 

We fix 
the measure $dz$ on $Z_K$ with total mass $1$, and we obtain a measure $dk$ on 
$K$ by requiring that $dk/dz$ be the measure on the compact open subgroup $K/Z_K\subset SO(V)(\A_f)$ induced by $dh_f$. 
This provides a normalization of the Haar measure on $H(\A_f)$ and hence a measure $dh'$ on $Z(\R)\back H(\A)$. 
Continuing the calculation above, and noting that $Z(\A)=Z(\Q)Z(\R)Z_K$, we have
$$\align
I(g';\tilde{\ph}_{KM}\tt\ph;V) &
=\frac12 \int_{H(\Q)Z(\A)\back H(\A)} \theta(g',h';\tilde\ph_{KM}\tt\ph)\,dh\\
\nass
\nass
{}&=\frac12 \int_{H(\Q)Z(\R)\back H(\A)} \theta(g',h;\tilde\ph_{KM}\tt\ph)\,dh'\tag4.45\\
\nass
\nass
{}&=\frac12 \sum_j\int_{H(\Q)Z(\R)\back H(\Q)H(\R)^+ h_jKh_j^{-1}} \theta(g',hh_j;\tilde\ph_{KM}\tt\ph)\,dh'\\
\nass
\nass
{}&=\frac14\, \vol(K) \sum_j\int_{\Gamma_jZ(\R)\back H(\R)^+} \theta(g',h_\infty h_j;\tilde\ph_{KM}\tt\ph)\,dh_\infty.
\endalign
$$
Here we have used the fact that $\ph$ is $K$--invariant. The extra factor of $\frac12$ in the last step arises 
from the fact that $Z(\Q)\cap K\simeq \{\pm1\}$. Finally, we 
normalize the measure $dh_\infty$ on $Z(\R)\back H(\R) = SO(V)(\R)$ by requiring that for $\phi\in C_c(D)$, 
$$\int_{Z(\R)\back H(\R)} \phi(h_\infty z_0) \, dh_\infty = (-1)^n\int_D \phi\cdot \O^n,\tag4.46$$
where $z_0\in D$ is the base point used in (4.43). 
Then, using (4.39),  we have
$$\align
I(g';\tilde{\ph}_{KM}\tt\ph;V) &=
(-1)^n\,\frac14\, \vol(K) \sum_j\int_{\Gamma_j\back D^+} \theta(g', h_j;\tilde\ph_{KM}\tt\ph)\,\O^n\\
\nass
\nass
{}&=(-1)^n\,\frac14\, \vol(K) \sum_j\int_{\Gamma_j\back D^+} \theta(g', h_j;\ph_{KM}\tt\ph)\wedge\O^{n-1}\tag4.47\\
\nass
\nass
{}&=(-1)^n\,\frac14\,\vol(K) \int_{X_K}\theta(g',\ph_{KM}\tt\ph)\wedge \O^{n-1}.
\endalign
$$
\qed\enddemo

{\bf Remark 4.18.} The same unfolding argument yields
$$\align
1 &=\int_{O(V)(\Q)\back O(V)(\A)} dh\\
\nass
\nass
{}& = (-1)^n\,\frac14\,\vol(K)\,\sum_j\int_{\Gamma_j\back D^+} \O^n\tag4.48\\
\nass
\nass
{}&=(-1)^n \,\frac14\,\vol(K)\,\vol(X_K, \O^n),
\endalign
$$
and thus the useful formula
$$\vol(K) = (-1)^n\,\frac{4}{\vol(X_K,\O^n)}.\tag4.49$$
The sign in (4.44) has been introduced to make $\vol(K)$ positive. 

Viewed as a differential form on $D\times H(\A_f)/K$, 
the theta form is given by 
$$\theta(g',h;\ph_{KM}\tt\ph) = \sum_{m\in \Q} \sum_{\scr x\in V(\Q) \atop\scr Q(x)=m}
\o(g')\ph_{KM}(x)\,\ph(h^{-1}x),\tag4.50$$ 
on the set $D\times h K$. 
Here note that $\o(g')\ph_{KM}(x)$ is a $(1,1)$--form on $D$ and that (4.48) is, in fact, the Fourier 
expansion of the theta form as a function on $G'_\A$. 
Let $\theta_m(g';\ph_{KM}\tt\ph)$ be the $m$--th Fourier coefficient, i.e., the partial sum over $x\in V(\Q)$ with $Q(x)=m$, and 
note that, since this form is itself $H(\Q)$--invariant, it defines a $(1,1)$--form on $X_K$. 

Recall from section 1 that, for $m\in \Q_{>0}$ and for any $\ph\in S(V(\A_f))^K$, there is a divisor $Z(m,\ph)=Z(m,\ph;K)$ on $X_K$, 
cf. (1.48) and (1.52). 
Also recall that the line bundle $\Cal L_D^\vee$ on $D$ descends to a line bundle $\Cal L^\vee$ on $X_K$.
\proclaim{Definition 4.19} The $\Cal L^\vee$--degree of a cycle $Z$ of codimension $r$ in $X_K$ is
$$\deg_{\Cal L^\vee}(Z) := \int_{Z} \O^{n-r},$$
where $\O$ is the first Chern form of $\Cal L^\vee$, as in Proposition~4.11 and Corollary~4.12.
\endproclaim

Note that if $X_K$ were compact and smooth, this would be simply $c_1(\Cal L^\vee)^{n-r}[Z]$, for the 
first Chern class $c_1(\Cal L^\vee)$. 

Also observe that, for $Z$ an irreducible subvariety,
$$(-1)^{n-r}\deg_{\Cal L^\vee}(Z) >0\tag4.51$$
since $-\O$ is a K\"ahler form on $X_K$. 

The following result is a consequence of the Thom form property of $\ph_{KM}$ (Theorem~4.1 of \cite{\kmillsonII}, 
and Theorem~2.1 of \cite{\kmcana}). As before, take $\tau=u+iv\in \frak H$, and write $q^m = e(m\tau)$. 
\proclaim{Theorem 4.20} For $m>0$, and for $g'_\tau\in G'_\R$, 
$$
\int_{X_K} \theta_m(g'_\tau,\ph_{KM}\tt\ph)\wedge \O^{n-1} = v^{(\ell+2)/2}\,\deg_{\Cal L^\vee}(Z(m,\ph;K))\cdot q^m,
$$
where $\ell=\frac{n}2-1$.
\endproclaim

{\bf Remark 4.21:} A key point here is that the cycle $Z(m,\ph)$ always has finite volume and the invariant form $\O^{n-1}$ 
is, in particular, bounded. Thus, Theorem~2.1 of \cite{\kmcana} can be applied, even when $X_K$ is non-compact. 
Alternatively, it is easy to obtain Theorem~4.20 by a direct calculation, using the integral formulas 
for the affine symmetric spaces, as used in the estimates in section~3 above. 

We now turn to the theta integral for the space $V'$. 

We fix a compact open subgroup $K'\subset O(V')(\A_f)$ such that 
$\ph'\in S(V'(\A_f))^{K'}$, and write
$$O(V')(\A) = \coprod_{j} O(V')(\Q)\, O(V')(\R)\, h_j K'.\tag4.52$$
Note that, since $V'$ is positive definite, the group
$$\Gamma_j= O(V')(\Q) \cap \bigg( O(V')(\R)h_j K' h_j^{-1}\bigg)\tag4.53$$
is finite; we set $e_j=|\Gamma_j|$. 
Again, we have a standard calculation, where we note that the Gaussian $\ph'_0$ is invariant under $O(V')(\R)$:
$$\align
I(g',\ph'_0\tt\ph';V') &=\int_{O(V')(\Q)\back O(V')(\A)} \theta(g',h;\ph'_0\tt\ph')\,dh\\
\nass
\nass
{}&= \sum_j \int_{\Gamma_j\back O(V')(\R) h_j K'}  \theta(g',h;\ph'_0\tt\ph')\,dh\tag4.54\\
\nass
\nass
{}&= \vol(O(V')(\R)K')\,\sum_j e_j^{-1} \theta(g',h_j;\ph'_0\tt\ph')\,dh.
\endalign
$$
If we take $g'=g'_\tau$, then
$$\o(g'_\tau)\ph'_0(x) = v^{(\ell+2)/2}\, e(Q'(x)\tau).\tag4.55$$ 
Note that 
$$\align
1&=\vol(O(V')(\Q)\back O(V')(\A), dh)\tag4.56\\
\nass
{}&= \vol(O(V')(\R)K')\,\sum_j e_j^{-1},
\endalign
$$
we have
$$\vol(O(V')(\R)K') = \bigg(\sum_j e_j^{-1}\bigg)^{-1}:=\mu(K'),\tag4.57$$
the mass of the $K'$--genus.  Thus we obtain the classical expression
$$I(g'_\tau,\ph'_0\tt\ph';V') = v^{(\ell+2)/2}\,\mu(K')\, \sum_j e_j^{-1} \theta(g',h_j;\ph'_0\tt\ph').\tag4.58$$
\proclaim{Proposition 4.22} For $m\in \Q$, let $q^m= e(m\tau)$, and recall that $\ell =\frac{n}2-1$. 
Then the  Fourier expansion of $I(g'_\tau,\ph'_0\tt\ph';V')$ is given by
$$I(g'_\tau,\ph'_0\tt\ph';V') =v^{(\ell+2)/2}\, \sum_{m\ge 0} r_{\ph'}(m)\, q^m,$$
where
$$r_{\ph'}(m) = \mu(K')\,\sum_j e_j^{-1}\,\bigg(\sum_{\scr x\in V'(\Q)\atop \scr Q'(x)=m}\ph'(h_j^{-1}x)\bigg).$$
In particular,  the constant term is 
$$v^{(\ell+2)/2}\,\ph'(0) = v^{(\ell+2)/2}\,\ph(0),$$
via matching. 
\endproclaim

The matching identity now amounts to:
\proclaim{Theorem 4.23} For $\ph\in S(V(\A_f))$ and $\ph'\in S(V'(\A_f))$ matching, 
and for the corresponding standard section $\P(s)$, with $\P(s_0) = \P_\infty^{\ell+2}(s_0)\tt\lambda(\ph)$, 
$$\align 
v^{-(\ell+2)/2}\,E(g'_\tau,s_0,\P) &= \ph(0) + \frac{1}{\vol(X_K,\O^n)} \sum_{m>0} \deg_{\Cal L^\vee}(Z(m,\ph;K))\cdot q^m\\
\nass
\nass
{}&=\ph'(0) + \sum_{m>0} r_{\ph'}(m)\, q^m.
\endalign
$$
\endproclaim

Comparing coefficients, we have
\proclaim{Corollary 4.24}
(i) For  $m>0$, 
$$\deg_{\Cal L^\vee}(Z(m,\ph;K)) = \vol(X_K,\O^n) \, r_{\ph'}(m).$$
(ii) If $m=0$, 
$$ \int_{X_K}\theta_0(g'_\tau,\ph_{KM}\tt\ph)\wedge \O^{n-1} =  v^{(\ell+2)/2}\, \ph(0)\,\vol(X_K,\O^n).$$
(iii) If $m<0$, 
$$\int_{X_K}\theta_m(g'_\tau,\ph_{KM}\tt\ph) = 0.$$
\endproclaim

Several special cases of these identities occur in the literature, cf. for example, \cite{\vdgeerbook}, \cite{\hermanna}, \cite{\hermannb}.
Note that, whereas the Fourier coefficients of the theta integral for $V$ involve 
degrees $\deg_{\Cal L^\vee}(Z(m,\ph))$, the Fourier coefficients of the theta integral 
for the positive definite space $V'$ are weighted representation numbers and the Fourier 
coefficients of the special value of the Eisenstein series, at least for factorizable data, have a product formula, i.e., 
are `multiplicative functions' in classical terminology.

\subheading{\Sec5. Examples}

In this section, we illustrate our results about integrals of Borcherds forms 
and about generating functions for degrees with some more or less explicit examples.
The basic idea is the following. On the one hand, by using Hasse--Minkowski, one can construct
even integral quadratic lattices $M$ of signature $(n,2)$ with prescribed local behavior. 
Associated to such a lattice are global geometric objects, the quasi--projective 
variety $X_M = \Gamma_M\back D^+$, the divisors
$Z(m,\ph)$, $\ph\in M^\sh/M$, etc.  On the other hand, associated to cosets
$\ph\in M^\sh/M$ are the Eisenstein series $E(\tau,s;\ph)$ of weight $\frac{n}2+1$.
These series and their Fourier expansions depend directly on the local data defining $M$.
The local and global objects are then related by the degree identity of Theorem~4.23,
$$\vol(X_M)\cdot E(\tau,\frac{n}2;\ph) = \vol(X_M)\cdot \ph(0) + \sum_{m>0} \deg_{\Cal L^\vee}(Z(m,\ph))\cdot q^m,$$
giving the first term of the Laurent expansion at $s=\frac{n}2$, and by 
Theorem~2.12, expressing the (log--norm) integrals of all Borcherds forms $\Psi(F)^2$ 
for $\C[L^\sh/L]$--valued $F$'s of weight $1-\frac{n}2$ in terms of the $\kappa_\ph(m)$'s arising from the 
second term of $E(\tau,s;\ph)$'s at $s=\frac{n}2$. 

A more systematic discussion of examples will be given in a sequel with Tonghai Yang, \cite{\ky}. 

First, we recall an example due to Gritsenko and Nikulin \cite{\GN}. Let
$$Q = \pmatrix {}&1&{}&{}&{}\\1&{}&{}&{}&{}\\{}&{}&2&{}&{}\\{}&{}&{}&{}&1\\{}&{}&{}&1&{}\endpmatrix,\tag5.1$$
and let $M = \Z^5$ with quadratic form, of signature $(3,2)$, defined by 
$$Q(x) =\frac12 {}^tx Q x.\tag5.2$$
The dual lattice of $M$ is $M^\sh=Q^{-1}M$ and $|M^\sh/M|=2$. 
Note that, if $x\in M^\sh$ with $Q(x) = m$, then 
$m\in \frac14 \Z$ and $4m \equiv 0,1\mod(4)$, depending on the $M$ coset $x+M$. 

As explained in \cite{\GN}, pp186--188, and \cite{\vdgeer}, there are compatible isomorphisms
$$\text{\rm GSpin}(M_{\R})\simeq \text{\rm Sp}_4(\R), \qquad D^+ \simeq \frak H_2,\tag5.3$$
such that 
$$\Gamma = SO^+(M) \simeq \text{\rm Sp}_4(\Z)/\{\pm1_4\}.\tag5.4$$
Let 
$$X= \Gamma\back D^+ \simeq \text{\rm Sp}_4(\Z)\back \H_2. \tag5.5$$
Recall that, in the tube domain model, our invariant form $\O=\ph_{KM}(0)$ 
is given by
$$\O = -\frac{1}{4\pi i} \bigg[\ -2(y,y)^{-2}(y,dz)\wedge (y,d\bar z) + (y,y)^{-1}\,(dz,d\bar z)\ \bigg].\tag5.6$$
In the case $n=3$, we write
$$z = \pmatrix z_1&z_2\\z_2&z_3\endpmatrix\in \H_2\tag5.7$$
and take
the inner product of a pair of $2\times 2$ symmetric matrices to be 
$(a,b) = -\tr(ab^\iota)$, for $\iota$ the main involution on $M_2(\Q)$. By an easy computation, noting that 
$(y,y)= -2\det(y)$, we find:
$$\O^3   = -\frac{3}{16\pi^3}\,\det(y)^{-3}\,\left(\frac{i}2\right)^3\,
dz_1\wedge d\bar z_1\wedge dz_2\wedge d\bar z_2\wedge dz_3\wedge d\bar z_3\tag5.8$$
and so \cite{\vdgeer}, p.331
$$\vol(X,\O^3) = \zeta(-1)\,\zeta(-3) = -\frac1{12}\,\zeta(-3)= -\frac1{1440}.\tag5.9$$

Let $V= M\tt_\Z\Q$ be the associated rational quadratic space, and let
$\ph_0, \ \ph_1\in S(V(\A_f))$ be the characteristic functions of the sets
$\hat M = M\tt\hat\Z$ and $y_1+ \hat M$ respectively, where $y_1$ is an 
element in $M^\sh\setminus M$. As explained in \cite{\GN} and \cite{\vdgeer}, the 
divisors $Z(m,\ph_\mu)$, for $\mu=0$, $1$,  are then given by
$$Z(m,\ph_\mu) = \cases \Cal G_{4m}&\text{ if $4m\equiv \mu\mod(4)$,}\\
\nass
\emptyset&\text{otherwise,}\endcases\tag5.10
$$
where 
$$\Cal G_{\Delta} = \sum_{n\ge1\atop n^2\mid \Delta} \nu_{\Delta/n^2}\,\Cal H_{\Delta/n^2},\tag5.11$$
for $\Cal H_\Delta$ the Humbert surface of discriminant $\Delta$ and with 
$$\nu_\Delta = \cases \frac12&\text{ if $\Delta=1$ or $4$,}\\
\nass
1&\text{ otherwise.}
\endcases\tag5.12
$$

We can define a vector valued Eisenstein 
series
$$\bold E(\tau,s;M) = \pmatrix E(\tau,s;\ph_0)\\ E(\tau,s;\ph_1)\endpmatrix\tag5.13$$
of weight $5/2$. The Fourier expansion of this series can be computed \cite{\ky}, and
from this it is easy to derive the following information. Write
$$E(\tau,s;\ph_\mu) = \sum_m A_\mu(s,m,v)\, q^m\tag5.14$$ 
as in (2.21), 
where the Fourier coefficients have Laurent expansions 
$$A_\mu(s,m,v) = a_\mu(m) + b_\mu(m,v)(s-s_0) + O((s-s_0)^2),\tag5.15$$
as in (2.22). 

\proclaim{Proposition 5.1} The value of $\bold E(\tau,\frac32;L_0)$ at the point $s_0=\frac32$ is given by the following expression.\hfill\break
$$\align
E(\tau,\frac32;\ph_0)  &= 1+ \zeta(-3)^{-1}\sum_{m=1}^\infty H(2,4m)\,q^m\\
\nass
\noalign{and}
\nass
E(\tau,\frac32;\ph_1) &= \phantom{1 + }\ \zeta(-3)^{-1}\sum_{m - \frac14 =0}^{\infty} H(2,4m)\,q^m
\endalign
$$
where $H(2,N)$ are as in Cohen \cite{\cohen}. \hfill\break
\endproclaim

In particular, for the value, observe that
$$E(4\tau,\frac32;\ph_0) + E(4\tau,\frac32;\ph_1) = \zeta(-3)^{-1}\Cal H_2(\tau),\tag5.16$$
is Cohen's Eisenstein series of weight $\frac52$.
Also, for convenient reference, we recall some values from \cite{\cohen}:
$$
\matrix
N: &0&1&4&5&8&9&12&13&16&17&20&\dots\\
\nass 
-120\,H(2,N): &-1&10&70&48&120&250&240&240&550&480&528&\dots
\endmatrix\tag5.17
$$

{\bf Remark 5.2.}  Recall that the positive coefficients in Cohen's Eisenstein series 
$\Cal H_r(\tau)$ of weight $r+\frac12$ are given by
$$H(r,4m) = L(1-r,\chi_d)\,\sum_{c\mid n} \mu(c)\,\chi_d(c)\, c^{r-1}\s_{2r-1}(n/c),\tag5.18$$
where $4m = (-1)^rn^2d$ for a field discriminant $d\equiv 0,\ 1\mod(4)$. The sum on $c$ 
is a multiplicative function and it is easy to check that, in fact, 
$$H(r,4m) = L(1-r,\chi_d)\,\prod_p b_p(n,1-r).\tag5.19$$
where
$b_p(n,s)$ is given by
$$b_p(n,s) = \frac{1-\chi_d(p) \,X +\chi_d(p)\,p^k X^{2k+1} - p^{k+1}X^{2k+2}}{1-p X^2},\tag5.20$$
with $X=p^{-s}$ and $k=\ord_p(n)$.

By Theorem~4.23, we have
$$
E(\tau,\frac32;\ph_\mu) =\ph_\mu(0)+ \vol(X)^{-1}\sum_{m>0} \deg(Z(m,\ph_\mu))\,q^m,\tag5.21
$$
so we obtain, for $4m\equiv \mu\mod(4)$, 
$$\deg(Z(m,\ph_\mu)) =\deg(\Cal G_{4m}) =  -\frac1{12}\,H(2,4m)\tag5.22$$
Thus, we recover the relation (1) of van der Geer, \cite{\vdgeer}, p.346, as well as his 
Theorem~8.1 on the generating function for the volumes of the Humbert surfaces. 

A nice example of a Borcherds form $\Psi(f)$ is discussed in \cite{\GN}.

Let $\phi_{12,1}(\tau,w)$, $\tau\in \H_1$, $w\in \C$ be the holomorphic Jacobi form of
weight $12$ and index $1$ of Eichler and Zagier \cite{\EZ}, pp.38--39, so that
$$\phi_{12,1}(\tau,w) = \sum_{n,r} C_{12}(4n-r^2)\,q^n\,\zeta^r,\tag5.23$$
for $q=e(\tau)$ and $\zeta=e(w)$, where $c_{12}(n)$ 
is given by the table on p.141 of \cite{\EZ}:
$$\matrix n&0&3&4&7&8&11&12&15&16&19&20&\dots\\
\nass
C_{12}(n)&0&1&10&-88&-132&1275&736&-8040&-2880&24035&13080&\dots\endmatrix\tag5.24
$$
(We write $C_{12}(n)$ in place of $c_{12}(n)$ to avoid confusion with the 
coefficients $c_\mu(m)$ which will occur in a moment.)
Write
$$\phi_{12,1}(\tau,w) = \sum_{\mu=\, 0,\, 1} h_\mu(\tau)\,\theta_{1,\mu}(\tau,w)\tag5.25$$
where
$$h_\mu(\tau) = \sum_{m\atop m\equiv -\mu\mod(4)} C_{12}(m)\,q^{\frac{m}4}\tag5.26$$
has weight $\frac{23}2$ for $\Gamma_0(4)$ and $\theta_{1,\mu}(\tau,w)$ is the standard Jacobi theta series. 
Then, dividing by $\Delta$ to shift the weight, we have
$$\frac{\phi_{12,1}(\tau,w)}{\Delta(\tau)} = \sum_{\mu=\, 0,\, 1} f_\mu(\tau)\,\theta_{1,\mu}(\tau,w)\tag5.27$$
where
$$f_\mu(\tau) = \sum_{m} c_\mu(m)\,q^m,\tag5.28$$
has weight $-\frac12$ and
$$\align
f_0(\tau) &= 10 + 108\, q + 808\, q^2 + ...\tag5.29\\
\nass
f_1(\tau) &= q^{-\frac14} - 64\, q^{\frac34} - 513\, q^{\frac74} + ...  
\endalign
$$
Associated to the vector valued form (see \cite{\borch}, Example 2.3, p.500 and 
\cite{\EZ}, Theorem~5.1, p.59) 
$$\bold f_{5}(\tau) = \big(\ f_0(\tau),f_1(\tau)\ ) = f_0(\tau)\,\ph_0+f_1(\tau)\,\ph_1,\tag5.30$$
valued in $\C[M^\sh/M]$, 
is a Borcherds form $\Psi(\bold f_{5})$, identified explicitly by Gritsenko-Nikulin:
$$\Psi(\bold f_5) = 2^{-6}\,\Delta_5(z),\tag5.31$$
where $\Delta_5(z)$ is the Siegel cusp form of weight $5$ (and character) for $\text{Sp}_4(\Z)$. 
Then $\Psi(\bold f_5)^2$ has weight $10$ (and trivial character) and 
$$\div(\Psi(\bold f_{5})^2) = Z(\frac14,\ph_1).\tag5.32$$

Similarly, for any positive integer $t$, we can consider the form
$j(\tau)^t\cdot  f(\tau)$.
For example, for $t=1$, we get
$$\align 
j(\tau)\,f_0(\tau) &= 10\,q^{-1}+7548+ O(q)\tag5.33\\
\nass 
j(\tau)\,f_1(\tau) & = q^{-\frac54}+680\,q^{-\frac14} + O(q^{\frac14})
\endalign
$$
so that the associated $\Psi(\bold f_{3774})^2$ has weight $7548$ and divisor
$$10\,Z(1,\ph_0) + Z(\frac54,\ph_1) + 680\,Z(\frac14,\ph_1).\tag5.34$$
For $t=2$, we get
$$\align 
j(\tau)^2\,f_0(\tau) &= 10\,q^{-2}+14988\,q^{-1} + 9634552+O(q)\tag5.35\\
\nass 
j(\tau)^2\,f_1(\tau) & = q^{-\frac94}+1424\,q^{-\frac54} + 851559\,q^{-\frac14} + O(q^{\frac14})
\endalign
$$
so that the associated $\Psi(\bold f_{4827376})^2$ has weight $9634552$ and divisor
$$10\,Z(2,\ph_0) + 14988\,Z(1,\ph_0) + Z(\frac94,\ph_1) + 1424\,\,Z(\frac54,\ph_1)+ 851559\,Z(\frac14,\ph_1).\tag5.36$$
It is amusing to check the weight/degree relation, (2.30), 
$$\sum_{\mu}\sum_{m>0} c_\mu(-m)\,\frac1{12}H(2,4m) = -\vol(X)\,c_0(0),\tag5.37$$
i.e., 
$$-\sum_\mu\sum_{m>0} c_\mu(-m)\,120\,H(2,4m) = c_0(0)\tag5.38$$
in these cases. 

To compute the quantities $\kappa(\Psi(f))$ for these Borcherds forms, we need to determine the 
quantities $\kappa_\mu(m)$ derived from the second term in the Laurent expansion of 
$\bold E(\tau,s;M)$ at the point $s=\frac32$.

\proclaim{Theorem 5.3} (i) For $m>0$, write $4m=n^2d$ for $d$ the discriminant of the real quadratic 
field $\Q(\sqrt{m})$, and let $\chi_d$ be the associated quadratic character\footnote{ When $4m=n^2$, we take $\Q(\sqrt{m}) = \Q\oplus \Q$, 
$\chi_d=1$ and $L(s,\chi_1)=\zeta(s)$}. Then, for $4m\equiv \mu\mod(4)$, 
$$\align
b_\mu(m,v) &=  \zeta(-3)^{-1}\,H(2,4m)\,\bigg[\frac43 \, 
+2\,\frac{\zeta'(-3)}{\zeta(-3)} -\frac12\log(d)-\frac{L'(-1,\chi_d)}{L(-1,\chi_d)} -C\\
\nass
\nass
{}&\qquad\qquad +\sum_{p\mid n} \bigg(\,\log|n|_p-\frac{b'_p(n,-1)}{b_p(n,-1)}\, \bigg) +\frac12\,J(\frac32,4\pi m v)\,\bigg]. 
\endalign
$$
where 
$$2C = \log(4\pi)+\gamma,$$
$$J(\frac32,t) = \int_0^\infty e^{-tr} \frac{(1+r)^{\frac32} -1}{r}\,dr,$$ 
\smallskip
and for $k=\ord_p(n)$, 
$$-\frac{1}{\log(p)}\cdot\frac{b'_p(n,-1)}{b_p(n,-1)} = \frac{2p^3}{1-p^3} + \frac{ -\chi_d(p)\,p + \chi_d(p)\,(2k+1)\,p^{3k+1} -(2k+2)\,p^{3k+3}}{
1-\chi_d(p)\,p + \chi_d(p)\,p^{3k+1} - p^{3k+3}}.$$
(ii) For $m<0$, 
$$b_\mu(m,v) = -\frac{\pi^2}3\,\frac{L(2,\chi_m)}{\zeta(4)}\,(\pi v)^{-\frac32}\int_1^\infty e^{-4\pi|m|vr}\,r^{-\frac32}\,dr.$$
(iii) For the constant term is given by
$$b_0(0,v) = \frac12\,\log(v) - \frac{\pi}6\,\frac{\zeta(3)}{\zeta(4)}\,v^{-\frac32}.$$
(iv) If $4m\not\equiv \mu\mod(4)$, then $b_\mu(m,v)=0$.
\endproclaim 
For $m<0$, the $L$--series $L(s,\chi_m)$ is a modified Dirichlet series analogous to that 
occurring in the definition of $H(r,4m)$. In any case, it is clear that, $\lim_{v\rightarrow\infty} b_\mu(m,v)=0$ 
for $m<0$. Similarly, for $m>0$, $\lim_{v\rightarrow\infty}J(\frac32,4\pi m v) =0$. 

\proclaim{Corollary 5.4} For $m>0$ with $4m=n^2d$ and with $4m\equiv \mu\mod(4)$, 
$$\align
\kappa_\mu(m) &=  \zeta(-3)^{-1}\,H(2,4m)\,\bigg[\frac43 \, 
+2\,\frac{\zeta'(-3)}{\zeta(-3)} -\frac12\log(d)-\frac{L'(-1,\chi_d)}{L(-1,\chi_d)} -C\\
\nass
\nass
{}&\qquad\qquad\qquad\qquad\qquad\qquad +\sum_{p\mid n} \bigg(\,\log|n|_p-\frac{b'_p(n,-1)}{b_p(n,-1)}\, \bigg) \,\bigg]. 
\endalign
$$
If $4m\not\equiv \mu\mod(4)$, then $\kappa_\mu(m)=0$. 
\endproclaim

Now, in calculating $\kappa(\Psi(\bold f))$ via Theorem~2.12, we can use the degree relation:
$$\align
\kappa(\Psi(\bold f)) &= \sum_{\mu}\sum_{m>0} c_\mu(-m)\,\kappa_\mu(m) + c_0(0)\,\frac12\,C_0\\
\nass
{}&= \sum_{\mu}\sum_{m>0} c_\mu(-m)\,120\,H(2,4m)\bigg[ -\frac12\log(d) -\frac{L'(-1,\chi_d)}{L(-1,\chi_d)}\tag5.39\\
\nass
\nass
{}&\qquad\qquad\qquad\qquad\qquad\qquad +\sum_{p\mid n} \bigg(\,\log|n|_p-\frac{b'_p(n,-1)}{b_p(n,-1)}\, \bigg) \,\bigg]\\
\nass
\nass
{}&\qquad\qquad - c_0(0)\,\bigg[ \frac43 \, 
+2\,\frac{\zeta'(-3)}{\zeta(-3)} -C -\frac12 C_0\,\bigg]. 
\endalign
$$

In the first example above, where $m=\frac14$, $d=1$, $\chi_d=1$ and $L(s,\chi_d)=\zeta(s)$, we obtain 
$$\align
\kappa(\Psi(\bold f_{5})) &= \zeta(-3)^{-1}\,H(2,1)\,\bigg[ \, \frac43 +2\,\frac{\zeta'(-3)}{\zeta(-3)}-\frac{\zeta'(-1)}{\zeta(-1)}
-C\,\bigg] + 10\cdot \frac12\,C_0.\tag5.40\\
\nass
{}&= 10\,\bigg[ - \frac43 -2\,\frac{\zeta'(-3)}{\zeta(-3)}+\frac{\zeta'(-1)}{\zeta(-1)}
+\frac32\,\log(2)+\log(\pi)\,\bigg].
\endalign
$$
Noting that $|y|^2 = 2\det(y)$ here, we have
$$||\Psi(\bold f_5)(z)||^2 = 2^{-12}\,|\Delta_5(z)|^2\,2^5\,\det(y)^{5},\tag5.41$$
so that 
$$-\vol(X)^{-1}\int_X \log\big(|\Delta_5(z)|^2\,\det(y)^{5})\cdot\O^3 = 10\,\bigg[ - \frac43 -2\,\frac{\zeta'(-3)}{\zeta(-3)}+\frac{\zeta'(-1)}{\zeta(-1)}
+\frac32\,\log(2)+\log(\pi)\,\bigg] - 7\log(2).\tag5.42$$

In the second example, there are terms for $m=\frac54$, $1$ and $\frac14$, and we obtain
$$
\align
\kappa(\Psi(\bold f_{3774}))  =\ &700\,\bigg[\, \frac{\zeta'(-1)}{\zeta(-1)} + \frac{b'_2(2,-1)}{b_2(2,-1)} + \log(2)\,\bigg] \\
\nass
{}&+ 48\,\bigg[\,\frac{L'(-1,\chi_5)}{L(-1,\chi_5)} + \frac12\,\log(5)\,\bigg]\tag5.43\\
\nass
{}&+ 6800\,\frac{\zeta'(-1)}{\zeta(-1)}\\
\nass
{}&+ 7548 \,\bigg[-\frac43 - 2\,\frac{\zeta'(-3)}{\zeta(-3)}
+\frac32\,\log(2)+\log(\pi)\,\bigg].
\endalign
$$
where
$$\frac{b'_2(2,-1)}{b_2(2,-1)} = -\frac{9}{11}\,\log(2).\tag5.44$$
And so on. 

In the next section, we explain why the values $\kappa_\mu(m)$ which occur here should be connected with the 
`arithmetic volumes' of (suitable integral extensions of) the cycles $Z(m,\ph_\mu)$. 

\subheading{\Sec6. Speculations}

The integrals considered in this paper play a role in the arithmetic geometry 
of cycles on the $\GSpin(n,2)$ varieties discussed above. While these Shimura varieties 
have canonical models over $\Q$, for all $n$, we do not have 
a sufficient theory of the integral models to give a precise discussion of the 
integral extensions of the $Z(m,\ph)$'s for general $n$. In addition, even for the 
archimedean theory, due to the non-compactness 
of $X_K$, one will need a suitable 
theory of line bundles with singular metrics, Green's currents with additional singularities, etc. 
Such problems are under consideration by Burgos, Kramer and K\"uhn
\cite{\burgoskuehn}. For the case of arithmetic surfaces,  i.e., $n=1$, see \cite{\bost}, \cite{\kuehn}.  
Nonetheless, based on low dimensional 
calculations, it is possible to make some rough speculations, which provide a setting 
for the results of this paper. 

A metrized line bundle $\hat\o$ on a projective arithmetic variety $\XX$ over $\Spec(\Z)$ 
defines a class $\hat\o\in\Pich(\XX)\simeq \CH^1(\XX)$ and classes $\hat\o^r\in \CH^r(\XX)$, 
the $r$-th arithmetic Chow group of $\XX$,with rational coefficients \cite{\gilletsoule}. 
For a cycle $\ZZ$ on $\XX$ of codimension $r$, there is a height $h_{\hat\o}(\ZZ)$ 
with respect to $\hat\o$, \cite{\bostgilletsoule}. For example, for an integral horizontal 
$\ZZ$ of codimension $r$, with normalization $j:\tilde\ZZ\rightarrow \ZZ\subset \bXX$, assumed to be 
itself regular over $\Spec(\Z)$, 
$$h_{\hat\o}(\ZZ) = \degh\, j^*(\hat\o^{n+1-r}),\tag6.1$$
where $\degh:\CH^{n+1-r}(\ZZ)\rightarrow\R$ is the arithmetic degree map.
Also, if $(\ZZ,g)\in \CH^r(\XX)$ is a codimension $r$ cycle with Green's current $g$, 
then, for the height pairing $\langle\ ,\ \rangle$ between $\CH^r(\XX)$ and $\CH^{n+1-r}(\XX)$, 
$$\langle (\ZZ,g), \hat\o^{n+1-r}\rangle = h_{\hat\o}(\ZZ) + \frac12\,\int_{\XX(\C)}g\cdot c_1(\hat\o)^{n+1-r},\tag6.2$$
where $c_1(\hat\o)$ is the first Chern form of $\hat\o$ on $\XX(\C)$. 

For $V$ of signature $(n,2)$, let $X= X_K$ be the canonical model over $\Q$ of the 
arithmetic quotient $\Gamma_K\back D^+$. Here we are assuming that $K$ is large enough so that 
$X$ is geometrically irreducible. Suppose that we have a regular model 
$\XX$ of $X$ over $\Spec(\Z)$, with a regular compactification $\bXX$. Suppose that
the metrized line bundle $\Cal L^\vee$ dual to $\Cal L$ (cf. (1.4) and (1.5)) on $X$ is the restriction 
of a line bundle $\hat\o$ on $\bXX$, where the metric on $\hat\o$ is allowed to have 
singularities along $\bar X(\C)\setminus X(\C)$. Note that the first Chern form of $\hat\o$ 
is the form $\O$ considered above. Suppose that one has a sufficiently extended theory 
of an arithmetic Chow ring (with rational coefficients) $\CH^\bullet(\bXX)$ 
so that the height construction can be applied. Thus, in particular, 
$\hat\o$ defines a class in $\CH^1(\bXX)$ and powers $\hat\o^r\in \CH^r(\bXX)$, etc. 

Next, we consider a Borcherds form $\Psi=\Psi(f)^2$ of weight $c_0(0)$. Then $\Psi$ is 
meromorphic function on $X(\C)\simeq \XX(\C)$, whose divisor is rational over $\Q$. We suppose that, 
in fact, there is a (rational) section $\tilde\Psi$ of $(\o^{-1})^{\tt c_0(0)}$ whose restriction to $\XX(\C)\simeq X(\C)$ 
is $\Psi$. It follows that $\divh(\tilde\Psi) = -c_0(0)\,\hat\o\in \CH^1(\bar\XX)$. 
Then, we would have
$$\align
-c_0(0)\,\langle\hat\o,\hat\o^n\rangle &= \langle\divh(\tilde\Psi),\hat\o^n\rangle\\
\nass
{}&= h_{\hat\o}(\div(\tilde\Psi)) + \frac12 \int_{X(\C)}\log||\Psi||^{-2}\,\O^n\tag6.3\\
\nass
{}&= h_{\hat\o}(\div(\tilde\Psi)) + \frac12 \vol(X)\,\kappa(\Psi).
\endalign
$$

Recall that (Theorem~1.3), on $X(\C)$, 
$$\div_{X}(\Psi)=\div_{\XX_\Q}(\Psi(f)^2) = \sum_\ph\sum_{m>0} c_\ph(-m)\,Z(m,\ph).\tag6.4$$
Then, on the integral model, we would have
$$\div_{\XX}(\tilde\Psi) = \sum_\ph\sum_{m>0} c_\ph(-m)\,\ZZ(m,\ph)+\big(\,\text{\rm vertical components}\,\big),\tag6.5$$
where the $\ZZ(m,\ph)$'s have generic fibers $\ZZ(m,\ph)_\Q=Z(m,\ph)$. 

Using the expression in Theorem~2.12 for $\kappa(\Psi)= 2\kappa(\Psi(f))$, we obtain
$$\align
-c_0(0)\,\langle\hat\o,\hat\o^n\rangle &= \sum_{\ph}\sum_{m>0} c_\ph(-m)\bigg[ h_{\hat\o}(\ZZ(m,\ph)) +\vol(X)\,\kappa_\ph(m)\bigg]\tag6.6\\
\nass
\nass
&{}\qquad+ \vol(X)\,c_0(0)\,\kappa_0(0) + 
\text{\rm contributions of vertical components}.
\endalign
$$
This (hypothetical) relation is suggestive. For example, if $c_0(0)=0$ so that $\divh(\tilde\Psi)=0$, 
we obtain 
$$0= \sum_{\ph}\sum_{m>0} c_\ph(-m)\bigg[ h_{\hat\o}(\ZZ(m,\ph)) +\vol(X)\,\kappa_\ph(m)\bigg] + 
\text{\rm contributions of vertical components},\tag6.7
$$
which suggests a close relation between $\kappa_\ph(m)$ and the height $h_{\hat\o}(\ZZ(m,\ph))$. 

In our example for $n=5$ from section 5, we can write
$$\align
\kappa_\mu(m)& = \vol(X)^{-1}\,\deg(Z(m,\ph))\,
\bigg[\,-\frac12\log(d)-\frac{L'(-1,\chi_d)}{L(-1,\chi_d)}+\sum_{p\mid n} \bigg(\,\log|n|_p-\frac{b'_p(n,-1)}{b_p(n,-1)}\, \bigg) \,\bigg]\tag6.8\\
\nass
\nass
{}&\qquad\qquad+\vol(X)^{-1}\,\deg(Z(m,\ph))\,
\bigg[\,\frac43 \, 
+2\,\frac{\zeta'(-3)}{\zeta(-3)}  -C\,\bigg]
\endalign
$$
so that (6.6) can be written as
$$\align
-c_0(0)\,\langle\hat\o,\hat\o^3\rangle &= \sum_{\ph}\sum_{m>0} c_\ph(-m)\,\delta(m,\ph_\mu)\\
\nass
\nass
&{}\qquad+ \vol(X)\,c_0(0)\,\bigg[\,\kappa_0(0) - \frac43 \, 
-2\,\frac{\zeta'(-3)}{\zeta(-3)}  +C\,\bigg]\tag6.9\\
\nass
\nass
&{}\qquad\qquad\qquad+\text{\rm contributions of vertical components}.
\endalign
$$
where
$$\delta(m,\ph_\mu) = h_{\hat\o}(\ZZ(m,\ph_\mu)) +
\deg(Z(m,\ph))\,
\bigg[\,-\frac12\log(d)-\frac{L'(-1,\chi_d)}{L(-1,\chi_d)}+\sum_{p\mid n} \bigg(\,\log|n|_p-\frac{b'_p(n,-1)}{b_p(n,-1)}\, \bigg) \,\bigg].\tag6.10
$$
Again, this suggests that
$$h_{\hat\o}(\ZZ(m,\ph_\mu)) \equiv 
-\deg(Z(m,\ph))\,
\bigg[\,-\frac12\log(d)-\frac{L'(-1,\chi_d)}{L(-1,\chi_d)}+\sum_{p\mid n} \bigg(\,\log|n|_p-\frac{b'_p(n,-1)}{b_p(n,-1)}\, \bigg) \,\bigg]\tag6.11$$
and
$$
\langle\hat\o,\hat\o^3\rangle \equiv 
\vol(X)\,\bigg[\, \frac43 \, 
+2\,\frac{\zeta'(-3)}{\zeta(-3)} -\frac32\,\log(2)-\log(\pi) \,\bigg],\tag6.12$$
where, in both relations, we have still to account for a possible linear combination 
of $\log(p)$'s coming from vertical components. In addition, it is possible to 
shift a term of the form 
$$\vol(X)^{-1} \,\deg(Z(m,\ph))\cdot A,\tag6.13$$
where $A$ is a constant independent of $\mu$ and $m$, between the two terms in (6.8), so there is some further ambiguity. 
It seems reasonable to expect that $A$ is a multiple of $\zeta'(-1)/\zeta(-1)$. This would be consistent 
with recent results of Bruinier and K\"uhn for certain Hilbert modular varieties, \cite{\bruinierkuehn}, K\"uhn's thesis \cite{\kuehnthesis}, 
and conjectures of Maillot and Roessler, \cite{\maillotroessler}. 
Recall that $\vol(X) =
\zeta(-1)\,\zeta(-3)$. 

Of course, this discussion is too vague with 
respect to integral models, compactifications, an extended theory of arithmetic Chow rings, 
and vertical contributions. 
Nonetheless, it explains the motivation for 
considering the quantities $\kappa(\Psi(f))$ and $\kappa_\ph(m)$ and their 
possible applications.  

Finally, it is worthwhile to compare the formula for $\Psi(f)$ of Theorem~2.12 with the following result 
of Rohrlich\footnote{I am indebted to Ulf Kuhn for calling this beautiful paper to my attention.}, 
\cite{\rohrlich}. Suppose that $h$ is a meromorphic modular form of weight $k$ for a
Fuchsian group
$\Gamma\subset \text{\rm PSL}_2(\R)$  having a cusp of width $1$ at $\infty$ such that $h$ is nonzero at every cusp of $\Gamma$ and has 
constant term $1$ at $\infty$. Let 
$$\Cal E(z,s) =\sum_{\Gamma_\infty\back \Gamma} \Im(\gamma(z))^s, \qquad \Re(s)>1.\tag6.14$$
Write
$$\Cal E(z,s) = \frac{\vol(\Gamma\back \frak H)}{s-1} + g(z) + O(s-1).\tag6.15$$
Then, setting
$$r = \frac{2\pi}{\vol(\Gamma\back \frak H)},\tag6.16$$
we have 
$$g(z) = -\frac1{2\pi} \log(y^r |D(z)| e^{-r})\tag6.17$$
for a holomorphic modular form $D$ of weight $2r$ for $\Gamma$ (with possible multiplier). 
Then
$$\eqalign{
\kappa(h) :&= -\frac1{\vol(\Gamma\back \frak H)}\,\int_{\Gamma\back \frak H} \log\big|h(z) y^{\frac{k}2}\big|^2\, d\mu(z)\cr
\nass
\nass 
{}&= \sum_{z\in \Gamma\back \frak H} 2\,\ord_z(h)\, \log(y^r|D(z)|).\cr}
\tag6.18
$$ 

In the case of a Shimura curve case, i.e, for $V$ anisotropic of signature $(1,2)$ 
and a Borcherds form $\Psi(f)$, we have 
$$\ord_z(\Psi(f)^2) = 
\sum_\ph\sum_{m>0} c_\ph(-m) \,\ord_z\big( Z(m,\ph)\big).\tag6.19$$
This suggests the following question. 
Is there a  modular form $D^B$ on $X_K$ of weight $2r$ 
such that
$$\kappa_\ph(m) = \sum_{z\in Z(m,\ph)} \log\big(y^r|D^B(z)|e^{-\frac12rC_0}\big),\tag6.20$$
for 
$r=\frac{2\pi}{\vol(\Gamma_M\back \frak H)}$?

One would then have
$$\eqalign{
&\sum_z 2\ \ord_z(\Psi(f))\,\log\big(y^r|D^B(z)|\big)\cr
\nass
\nass
{}&=\sum_z\bigg(\sum_\ph\sum_{m>0} c_\ph(-m) \,\ord_z\big( Z(m,\ph)\big)\bigg)\,\log\big(y^r|D^B(z)|e^{-\frac12rC_0}\big)
+r C_0 \,\deg(\div(\Psi(f)))\cr
\nass
\nass
{}&=\sum_\ph\sum_{m>0} c_\ph(-m) \bigg(\sum_{z\in Z(m,\ph)}\log\big(y^r|D^B(z)|e^{-\frac12rC_0}\big)\ \bigg) +\frac12 C_0 c_0(0)\cr
\nass
\nass
{}&=\sum_\ph\sum_{m>0} c_\ph(-m)\,\kappa_\ph(m)  +\frac12 C_0 c_0(0)\cr
\nass
\nass
{}&=\sum_\ph\sum_{m\ge 0} c_\ph(-m)\,\kappa_\ph(m)\cr
\nass
\nass
{}&=\kappa(\Psi(f)),\cr
\nass
\nass
{}&=-\frac{1}{\vol(X_K)} \int_{X_K}\log||\Psi(z;f)||^2\,d\mu(z)
}\tag6.21
$$
just as in Rohrlich's case.
The function $D$ in Rohrlich arises in the Kronecker limit formula for the Fuchsian group $\Gamma$,  
\cite{\goldstein}. The function $D^B$ would be a kind of analogue in the compact quotient case, 
i.e., in the absence of the Eisenstein series!

\vskip .5in

\redefine\vol{\oldvol}

\bye